\documentclass[a4paper,12pt,onecolumn]{article}

\usepackage[T1]{fontenc}
\usepackage[utf8]{inputenc}
\usepackage[english]{babel}
\usepackage{amssymb,amsmath,amsthm,array,epsfig,lmodern,gensymb,pstricks,subfig,
hyperref,verbatim,booktabs,textcomp}

\newcommand{\R}{\mathbb R}
\newcommand{\D}{\mathcal D}
\newcommand{\W}{\vec{\mathcal W}}
\newcommand{\vol}{\mathcal{L}^d}
\newcommand{\surfvol}{\mathcal{H}^{d-1}}
\newcommand{\dH}[1]{\;{\rm d}{\cal H}^{#1}} 
\newcommand{\dL}[1]{\;{\rm d}{\cal L}^{#1}} 
\newcommand{\bigchi}{\ensuremath{\mathrm{\mathcal{X}}}}
\newcommand{\charfcn}[1]{\bigchi_{#1}} 
\newcommand{\Vh}{\underline{V}(\Gamma^m)}
\newcommand{\Wh}{W(\Gamma^m)}
\newcommand{\Vht}{\underline{V}(\Gamma^h(t))}
\newcommand{\Wht}{W(\Gamma^h(t))}
\newcommand{\uspacesimple}{\mathbb{U}}
\newcommand{\uspace}[1]{\mathbb{U}(\vec{#1})}
\newcommand{\uspacedisc}[2]{\mathbb{U}^{#2}(\vec{#1})}
\newcommand{\pspace}{\mathbb{P}}
\newcommand{\pnormspace}{\widetilde\pspace} 
\newcommand{\uspaceref}[1]{\mathcal{U}(\vec{#1})} 

\newcommand{\pspaceref}{\mathcal{P}} 
\newcommand{\pnormspaceref}{\widetilde\pspaceref}
\newcommand{\uspaceale}[1]{\mathfrak{U}(\vec{#1})} 
\newcommand{\uspacesemidiscale}[3]{\mathfrak{U}^{#2}(\vec{#1},#3)} 
\newcommand{\uspacediscale}[2]{\mathfrak{U}^{#2}(\vec{#1})} 
\newcommand{\pspaceale}{\mathfrak{P}} 
\newcommand{\pnormspaceale}{\widetilde\pspaceale}

\newcommand{\sigmaO}{o}
\newcommand{\nabs}{\nabla_{\!s}}
\newcommand{\id}{\rm id}
\newcommand{\ddt}{\frac{\rm d}{{\rm d}t}}
\newcommand{\NbulkT}{\vec{N}_{\Gamma,\Omega}^T}
\newcommand{\Nbulk}{\vec{N}_{\Gamma,\Omega}}
\newcommand{\errorXx}{\|\Gamma^h - \Gamma\|_{L^\infty}}
\newcommand{\LerrorUu}[1]{\|\vec U - I^h_{#1}\,\vec u\|_{L^2}}
\newcommand{\HerrorUu}[1]{\|\vec U - I^h_{#1}\,\vec u\|_{L^2(H^1)}}
\newcommand{\LerrorPp}{\|P - p\|_{L^2}}
\newcommand{\unitn}{\vec{\rm n}}
\newcommand{\unitt}{\vec{\rm t}}
\newcommand{\mat}[1]{\underline{\underline{#1}}\rule{0pt}{0pt}}
\newcommand{\V}{\vec{\mathcal{V}}} 
\newcommand{\strikes}{\mbox{$s\!\!\!\!\:/$}}
\newcommand{\myurl}[1]{\hfil\penalty 100 \hfilneg \hbox{\url{#1}}} 

\newcommand{\schemeAex}{{\mathcal{A}^h_{\rm ex}}}
\newcommand{\schemeAim}{{\mathcal{A}^h_{\rm im}}}
\newcommand{\schemeB}{{\mathcal{B}^h}}
\newcommand{\schemeALE}{{\mathcal{C}^h_{\rm ALE}}}
\newcommand{\ek}{{\rm e}}

\newtheorem{thm}{Theorem}

\newtheorem{rem}[thm]{Remark}

\newenvironment{keywords}{{\upshape\bfseries Key words. }\ignorespaces}{}

\title{Fitted front tracking methods for two-phase incompressible
Navier--Stokes flow: \\
Eulerian and ALE finite element discretizations
}
\author{Marco Agnese and Robert N\"urnberg%
\thanks{email: \texttt{\{m.agnese13,robert.nurnberg\}@imperial.ac.uk}}\\
\small
Department of Mathematics, Imperial College, London, SW7 2AZ, U.K.}
\date{}

\begin{document}

\captionsetup[subfigure]{labelformat=empty} 

\maketitle

\begin{abstract}
We investigate novel fitted finite element approximations for two-phase
Navier--Stokes flow. In particular, we consider both Eulerian and
Arbitrary Lagrangian--Eulerian (ALE) finite element formulations.
The moving interface is approximated with the help of
parametric piecewise linear
finite element functions. The bulk mesh is fitted to the interface
approximation, so that standard bulk finite element spaces can be used
throughout.
The meshes describing the discrete interface in
general do not deteriorate in time, which means that in numerical simulations a
smoothing or a remeshing of the interface mesh is not necessary.
We present several numerical experiments, including convergence experiments
and benchmark computations, for the introduced numerical methods,
which demonstrate the accuracy and robustness of the proposed algorithms.
We also compare the accuracy and efficiency of the Eulerian and ALE
formulations.
\end{abstract}


\section{Introduction}\label{sec:introdution}
Fluid flow problems with a moving interface are encountered in many
applications in physics, engineering and biophysics. Typical applications
include drops and bubbles, die swell, dam break, liquid storage tanks,
ink-jet printing and fuel injection. For this reason, developing robust and
efficient numerical methods for these flows is an important problem and
has attracted tremendous interest over the last few decades.

A crucial aspect of these types of fluid flow problems is that apart from
the solution of the flow in the bulk domain, the position of the interface
separating the two bulk phases also needs to be determined.
At the interface certain boundary conditions need to be fulfilled, which
specify the motion of the phase boundary.
These conditions relate the variables of the bulk flow, velocity and pressure,
across the two phase, taking into account external influences, such as for
example surface tension. Numerically, in order
to be able to compute the flow solution as well as the interface
geometry, a measure to track the interface starting from an initial position
needs to be incorporated.
There are several strategies to deal with this problem, which can be divided
into two categories depending on the viewpoint: interface capturing and
interface tracking.

Interface capturing methods use an Eulerian description of the interface, which
is defined implicitly. A characteristic scalar field,
that is advected by the flow, is used to identify
the two phases as well as the interface along the boundaries of the individual
fluid domains. Depending on the method, this scalar field may be, for example,
a discontinuous Heaviside function or a signed-distance function. The most
important methods, which belong to this category, are
the volume-of-fluid method, the level-set method and the phase-field method.
In the volume-of-fluid method, the characteristic
function of one of the phases is approximated numerically, see e.g.
\cite{HirtN81,RenardyR02,Popinet09}. In the level-set method, the interface is
given as the level set of a function, which has to be determined, see e.g.
\cite{SussmanSO94,Sethian99,OsherF03,GrossR07,GrossR11,Svacek17}. Instead, the
phase-field method works with diffuse interfaces, and therefore the transition
layer between the phases has a finite size. There is no tracking mechanism for
the interface, but the phase state is included implicitly in the governing
equations. The interface is associated with a smooth, but highly localized
variation of the so-called phase-field variable.
Examples for phase-field methods applied to two-phase flow are
\cite{HohenbergH77,AndersonMW98,LowengrubT98,Boyer02,Feng06,DingSS07,%
KaySW08,AbelsGG12,AlandV12,GrunK14,GarckeHK16}.
Extensions of the method to multi-phase flows
can be found in \cite{Dong14,Dong15,BoyerM14,BanasN17}.
The appeal of interface capturing approaches is the fact that they are usually
easy to implement, and that they offer an automated way to deal with
topological changes.

Interface tracking approaches, instead, use a Lagrangian description of the
interface which is described explicitly. Here the interface is represented
by a collection of particles or points, and this representation is
transported by the bulk flow velocity.
The great
advantage of interface tracking approaches is that they offer an accurate and
computationally efficient approximation of the evolving interface.
Challenges are the mesh quality both of the interface representation and of the
bulk triangulations, as well as the need to heuristically deal with topological
changes.
We refer e.g.\ to
\cite{UnverdiT92,Bansch01,Tryggvason_etal01,GanesanT08,spurious,
fluidfbp} for further details, and to \cite{LevequeL97,Peskin02} for the
related immersed boundary method, which is used to simulate fluid-structure
interactions using Eulerian coordinates for the fluid and Lagrangian
coordinates for the structure.

In this paper we use a direct description of the interface using a
parametrization of the unknown interface, similarly to the previous work
\cite{fluidfbp}. In particular, our numerical method will be based
on a piecewise linear parametric finite elements, and the description of the
interface will be advected in normal direction
with the normal part of the fluid velocity. The tangential degrees
of freedom of the interface velocity are implicitly used to ensure a
good mesh quality, and this is a main feature of our proposed methods.
But in contrast to \cite{fluidfbp}, where an unfitted bulk finite element
approximation was used, we will adopt a fitted mesh approach, which means
that the discrete interface is composed of faces of elements from the bulk
triangulation. Fitted and unfitted bulk mesh approaches are fundamentally
different approximation methods, and need specialized implementation
techniques.
A fitted method has the advantage that discontinuity jumps
at the interface are captured naturally, but it has the disadvantage that in
a standard Eulerian method the velocity needs to be interpolated from
an old mesh to a new mesh, unless the interface is stationary. Hence so-called
Arbitrary Lagrangian Eulerian (ALE) methods are often proposed. Here
the equations are posed in a moving domain framework,
and an arbitrary Lagrangian velocity may be chosen to improve the quality of
the bulk mesh. The original idea for ALE methods goes back to
the papers \cite{Donea83,Hughes81}, and we refer to
\cite{Nobile99,NobilePhd,Formaggia04,Ganesan06,GanesanT08,HahnHT13,%
GanesanHST17} for applications to two-phase Navier--Stokes flow.

In this paper we will propose both standard Eulerian and
ALE finite element approximations for two-phase incompressible Navier--Stokes
flow. One aim of our paper is to investigate numerically, if there is a clear
advantage of ALE type methods over standard Eulerian methods.
We stress that we are not aware of any detailed comparisons between fitted
Eulerian and ALE front tracking methods in the literature.
This paper aims to fill this gap. On the other hand, similar
comparisons between ALE type interface tracking methods and various
interface capturing methods have been presented in e.g.\
\cite{HysingTKPBGT09,ElgetiS16}.
We note that in the case of viscous incompressible two-phase Stokes flow,
there is no need for the numerical method to interpolate the velocity from the
old to the new mesh. Hence there is no need to employ an ALE method. In fact,
in the case of two-phase Stokes flow, all our proposed numerical methods
collapse the approximation considered by the authors in \cite{stokesfitted}.

The remainder of the paper is organized as follows. In
Section~\ref{sec:ns_model} we introduce the precise mathematical model of
two-phase flow that we study, and formulate the weak formulations on which our
Eulerian and ALE finite element approximations will be based. In
Section~\ref{sec:ns_fem} we introduce four different fully discrete finite
element discretizations, and in Section~\ref{sec:solution_method} we discuss
computational issues for their implementation as well as the employed solution
methods. Finally, in Section~\ref{sec:nr} we present several numerical
simulations for the introduced finite element approximations. In particular, we
investigate their accuracy in convergence tests with the help of two exact
solutions. Moreover, we compare the performance of our schemes in two well
known benchmark problems.

\section{Mathematical model and weak formulations}\label{sec:ns_model}
\begin{figure}
\centering
\includegraphics[width=.35\textwidth]{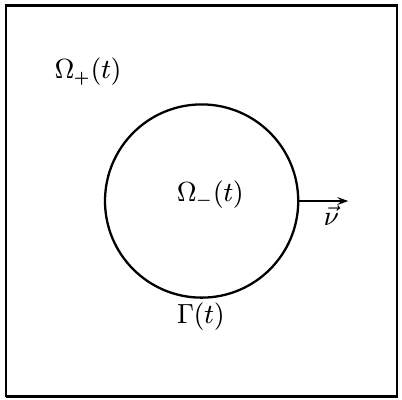}
\caption{The domain $\Omega$ in the case $d=2$.}
\label{fig:two_phase_sketch}
\end{figure}

Let $\Omega\subset\R^d$, for $d=2$ or $d=3$ be a given domain that contains two
different immiscible, incompressible fluids (liquid-liquid or liquid-gas) which
for all $t\in[0,T]$ occupy time dependent regions $\Omega_+(t)$ and
$\Omega_-(t):=\Omega\setminus\overline{\Omega}_+(t)$ and which are separated by
an interface $(\Gamma(t))_{t\in[0,T]}$, $\Gamma(t)\subset\Omega$. We limit
ourselves to interfaces formed by closed hypersurfaces, see
Figure~\ref{fig:two_phase_sketch} for a pictorial representation in the case
$d=2$.
For later use, we assume that $(\Gamma(t))_{t\in [0,T]}$ is a sufficiently
smooth evolving hypersurface without boundary that is parametrized by
$\vec x(\cdot,t):\Upsilon\to\R^d$, where $\Upsilon\subset \R^d$ is a given
reference manifold such that $\Gamma(t) = \vec x(\Upsilon,t)$. Then
\begin{equation} \label{eq:V}
\V(\vec z, t) := \vec x_t(\vec q, t) \quad
\forall\ \vec z = \vec x(\vec q,t) \in \Gamma(t)\,,
\end{equation}
defines the velocity of $\Gamma(t)$, and $\V \,.\,\vec \nu$ is
the normal velocity of the evolving hypersurface $\Gamma(t)$,
where $\vec\nu(t)$ is the unit normal on $\Gamma(t)$ pointing into
$\Omega_+(t)$.

Denoting the velocity and pressure by $\vec u$ and $p$, respectively, we
introduce the stress tensor
\begin{equation} \label{eq:stress_tensor}
\mat\sigma = \mu \,(\nabla\,\vec u + (\nabla\,\vec u)^T) - p\,\mat\id
= 2\,\mu\, \mat D(\vec u)-p\,\mat\id\,,
\end{equation}
where $\mu(t) = \mu_+\,\charfcn{\Omega_+(t)} + \mu_-\,\charfcn{\Omega_-(t)}$,
with $\mu_\pm \in \R_{>0}$, denotes the dynamic viscosities in the two phases,
$\mat\id \in \R^{d \times d}$ is the identity matrix and
$\mat D(\vec u):=\frac12\, (\nabla\,\vec u+(\nabla\,\vec u)^T)$
is the rate-of-deformation tensor. Here and throughout, $\charfcn{\D}$ denotes
the characteristic function for a set $\D$.

The fluid dynamics in the bulk domain $\Omega$ is governed by the two-phase
Navier--Stokes model
\begin{subequations}
\begin{alignat}{2}
\rho\,(\vec u_t +(\vec u \,.\, \nabla)\,\vec u)
- 2 \mu\,\nabla\,.\,\mat D(\vec u)+
\nabla\,p & = \vec f \quad &&\mbox{in } \Omega_\pm(t)\,,
\label{eq:ns_momentum_bis} \\
\nabla\,.\,\vec u & = 0 \quad &&\mbox{in } \Omega_\pm(t)\,,
\label{eq:ns_mass_bis}
\end{alignat}
\end{subequations}
where $\rho(t) = \rho_+\,\charfcn{\Omega_+(t)}+\rho_-\,\charfcn{\Omega_-(t)}$,
with $\rho_\pm \in \R_{>0}$, denotes the fluid density in the two phases, and
where $\vec f:=\rho\,\vec f_1+\vec f_2$ is a possible forcing term.
The velocity and stress tensor need to be coupled across the free surface
$\Gamma(t)$. Viscosity of the fluids leads to the continuity condition
\begin{equation}\label{eq:interface_jump_velocity}
[\vec u]_-^+ = \vec 0 \quad \mbox{on } \Gamma(t)\,,
\end{equation}
where $[\vec u]_-^+ := \vec u_+ - \vec u_-$ denote the jump in velocity across
the interface $\Gamma(t)$. Here and throughout, we employ the shorthand notation
$\vec h_\pm := \vec h\!\mid_{\Omega_\pm(t)}$ for a function
$\vec h : \Omega \times [0,T] \to \R^d$; and similarly for scalar and
matrix-valued functions. The assumption that there is no change of phase leads
to the dynamic interface condition
\begin{equation}\label{eq:interface_velocity}
\V\,.\,\vec\nu = \vec u\,.\,\vec \nu \quad \mbox{on }
\Gamma(t)\,,
\end{equation}
which means that the normal velocity of the interface needs to match the flow
normal velocity across the interface $\Gamma(t)$. Moreover, the momentum
conservation in a small material volume that intersects the interface leads to
the stress balance condition
\begin{equation}\label{eq:interface_jump_stress_div}
[\mat\sigma\,\vec \nu]_-^+ = -\nabs\,.\mat\sigma_{\Gamma} \quad \mbox{on }
\Gamma(t)\,,
\end{equation}
where $\mat\sigma_{\Gamma}$ is the interface stress tensor and $\nabs\,.$ is
the surface divergence on $\Gamma(t)$. The operator $\nabs$ is the surface
gradient on $\Gamma(t)$ and it is the orthogonal projection of the usual
gradient operator $\nabla$ on the tangent space of the surface $\Gamma(t)$.
Therefore, given a smooth function $h$ defined on a neighbourhood of
$\Gamma(t)$, it is defined as
\begin{equation}\label{eq:surface_gradient}
\nabs\,h(\vec z)=\mat P(\vec z)\,\nabla\,h(\vec z)=\nabla\,h(\vec z)\,-
\,\nabla\,h(\vec z)\,.\,\vec\nu(\vec z)\,\vec\nu(\vec z)\,,\quad \vec z
\in\Gamma(t)\,,
\end{equation}
with the usual projection operator $\mat P$ defined as
\begin{equation}\label{eq:surface_projection}
\mat P=\mat\id\,-\,\vec\nu\,\vec\nu\,^T\quad \mbox{on } \Gamma(t)\,.
\end{equation}
It can be shown that the definition of $\nabs\,h$ does not depend on the chosen
extension of $h$, but only on the value of $h$ on $\Gamma(t)$. See e.g.\
\cite[\S2.1]{DeckelnickDE05} for details.
For later use, we also define the surface
Laplacian, also known as Laplace--Beltrami operator, $\Delta_s$ on $\Gamma(t)$
as
\begin{equation}\label{eq:surface_laplacian}
\Delta_s = \nabs\,.\,\nabs\,.
\end{equation}
We restrict ourselves to the case that the only force acting on the interface
is the surface tension contact force. Therefore the interface stress tensor
$\mat\sigma_\Gamma$ has the following constitutive equation
\begin{equation}\label{eq:interface_stress_tensor}
\mat\sigma_{\Gamma}=\gamma\,\mat P\,,
\end{equation}
where $\gamma$ is the surface tension coefficient. Substituting
(\ref{eq:interface_stress_tensor}) in (\ref{eq:interface_jump_stress_div}) we
obtain, assuming $\gamma$ is constant, the so called clean interface model for
the interfacial forces
\begin{equation}
[\mat\sigma\,\vec \nu]_-^+ = -\gamma\,\varkappa\,\vec\nu \quad \mbox{on }
\Gamma(t)\,,
\end{equation}
which, in the case of Newtonian stress tensor (\ref{eq:stress_tensor}), assumes
the form
\begin{equation}\label{eq:interface_jump_stress}
[2\mu \,\mat D(\vec u)\,.\,\vec\nu - p\,\vec \nu]_-^+
= -\gamma\,\varkappa\,\vec\nu \quad \mbox{on } \Gamma(t)\,.
\end{equation}
Here $\varkappa$ denotes the mean curvature of $\Gamma(t)$, i.e.\ the sum of
the principal curvatures of $\Gamma(t)$, where we have adopted the sign
convention that $\varkappa$ is negative where $\Omega_-(t)$ is locally convex.
In particular, see e.g.\ \cite{DeckelnickDE05}, it holds that
\begin{equation} \label{eq:LBop}
\Delta_s\, \vec \id = \varkappa\, \vec\nu \qquad \mbox{on $\Gamma(t)$}\,.
\end{equation}

In order to close the system, we prescribe the
initial data $\Gamma(0) = \Gamma_0$ and the initial velocity
$\vec u(0) = \vec u_0$. Finally, we impose the Dirichlet condition $\vec u =
\vec g$ on $\partial_1 \Omega$ and the free-slip condition $\vec u \,.\,\unitn =
0$ and $\mat\sigma\,\unitn\,.\,\unitt = 0$, $\forall\ \unitt \in
\{\unitn\}^\perp$, on $\partial_2 \Omega$, with $\unitn$ denoting the outer
unit normal of $\partial \Omega$ and $\{\unitn\}^\perp := \{ \unitt \in \R^d :
\unitt \,.\,\unitn = 0\}$. As usual, it holds that
$\partial\Omega =\partial_1\Omega \cup \partial_2\Omega$ and $\partial_1\Omega
\cap \partial_2\Omega = \emptyset$. Therefore the total system of equations
can be rewritten as follows:
\begin{subequations}
\begin{alignat}{2}
\rho\,(\vec u_t +(\vec u \,.\, \nabla)\,\vec u)
- 2 \mu\,\nabla\,.\,\mat D(\vec u)+
\nabla\,p & = \vec f
\quad &&\mbox{in } \Omega_\pm(t)\,, \label{eq:ns_full_momentum} \\
\nabla\,.\,\vec u & = 0 \quad &&\mbox{in } \Omega_\pm(t)\,,
\label{eq:ns_full_mass} \\
\vec u = \vec g \quad \mbox{on } \partial_1\Omega\,,\qquad
\vec u \,.\,\unitn = 0 \,,\quad \mat\sigma\,\unitn\,.\,\unitt & = 0
\quad\forall\ \unitt \in \{\unitn\}^\perp \quad && \mbox{on } \partial_2\Omega\,,
\label{eq:ns_full_freeslip}\\
[\vec u]_-^+ & = \vec 0 \quad &&\mbox{on } \Gamma(t)\,,
\label{eq:ns_full_jump_velocity} \\
[2\mu \,\mat D(\vec u)\,.\,\vec\nu - p\,\vec \nu]_-^+
& = -\gamma\,\varkappa\,\vec\nu
\quad &&\mbox{on } \Gamma(t)\,, \label{eq:ns_full_jump_stress} \\
(\V-\vec u)\,.\,\vec \nu & = 0
\quad &&\mbox{on } \Gamma(t)\,,\label{eq:ns_full_velocity} \\
\Gamma(0) & = \Gamma_0 \,,\qquad \vec u(0) = \vec u_0 \,.
\end{alignat}
\end{subequations}
In what follows we introduce different weak formulations for
(\ref{eq:ns_full_momentum}--g), on which our finite element approximations will
be based.

\subsection{Eulerian weak formulations}\label{sec:ns_weak}
In order to obtain a weak formulation for (\ref{eq:ns_full_momentum}--g),
we define the function spaces, for a given $\vec b \in [H^1(\Omega)]^d$,
\begin{align*}
\uspace b &:= \{\vec \phi\in[H^1(\Omega)]^d:
\vec \phi =\vec b\quad \mbox{on }\partial_1\Omega\,,
\quad \vec \phi\,.\,\unitn=0 \quad \mbox{on }\partial_2\Omega\}\,,\\
\pspace &:= L^2(\Omega)\,,\qquad
\pnormspace := \{\eta\in\pspace : \int_\Omega\eta\dL{d}=0 \}\,,
\end{align*}
and let, as usual, $(\cdot,\cdot)$ and $\langle \cdot, \cdot
\rangle_{\Gamma(t)}$ denote the $L^2$--inner products on $\Omega$ and
$\Gamma(t)$, respectively. In addition, we let $\vol$ and $\surfvol$ denote the
Lebesgue measure in $\R^d$ and the $(d-1)$-dimensional Hausdorff measure,
respectively.

We also need a weak form of the differential geometry identity
(\ref{eq:LBop}), which can be obtained by multiplying (\ref{eq:LBop}) with a
test function and performing integration by parts leading to
\begin{equation}\label{eq:weak_nabs_id}
\left\langle \varkappa\,\vec\nu, \vec\eta \right\rangle_{\Gamma(t)}
+ \left\langle \nabs\,\vec \id, \nabs\,\vec \eta \right\rangle_{\Gamma(t)}
= 0 \quad \forall\ \vec\eta \in [H^1(\Gamma(t))]^d\,.
\end{equation}
Moreover, on noting (\ref{eq:stress_tensor}) and
(\ref{eq:interface_jump_stress}), we have that
\begin{align}\label{eq:weak_stress_jump}
& \int_{\Omega_+(t)\cup\Omega_-(t)} (\nabla\,.\,\mat\sigma)\,.\, \vec \xi \dL{d}
= - \left(\mat\sigma, \nabla\,\vec \xi\right)
- \left\langle [\mat\sigma\,\vec\nu]_-^+ , \vec \xi \right\rangle_{\Gamma(t)}
\nonumber \\
& \qquad = - \left(2\,\mu\, \mat D(\vec u)-p\,\mat\id, \nabla\,\vec \xi\right)
- \left\langle [2\mu \,\mat D(\vec u)\,.\,\vec\nu - p\,\vec \nu]_-^+ , \vec \xi
\right\rangle_{\Gamma(t)} \nonumber \\
& \qquad = \left( p , \nabla\,.\,\vec \xi\right)
-2 \left(\mu\,\mat D(\vec u) , \mat D(\vec \xi) \right)
+ \gamma \left\langle \varkappa\,\vec \nu , \vec \xi \right\rangle_{\Gamma(t)}
\end{align}
for all $\vec \xi \in \uspace 0$.
Hence a standard weak formulation of (\ref{eq:ns_full_momentum}--g)
is given as follows. Given $\Gamma(0) = \Gamma_0$ and
$\vec u(\cdot,0) = \vec u_0$, for
almost all $t\in(0,T)$ find $\Gamma(t)$ and $(\vec u, p, \varkappa)(t)\in
\uspace g \times \pnormspace \times H^1(\Gamma(t))$ such that
\begin{subequations}
\begin{align}
& \left(\rho\,\vec u_t, \vec \xi\right) + \left(\rho\,(\vec u \,.\, \nabla)\,
\vec u,\vec \xi\right) + 2\left(\mu\,\mat D(\vec u), \mat D(\vec \xi)\right)
- \left(p, \nabla\,.\,\vec \xi\right)
- \gamma\,\left\langle \varkappa\,\vec\nu, \vec\xi\right\rangle_{\Gamma(t)}
= \left(\vec f, \vec \xi\right)
\nonumber \\ & \hspace{11cm}
\forall\ \vec\xi \in \uspace 0 \,, \label{eq:ns_weaka}\\
& \left(\nabla\,.\,\vec u, \varphi\right) = 0
\quad \forall\ \varphi \in \pnormspace\,, \label{eq:ns_weakb} \\
& \left\langle \V
- \vec u, \chi\,\vec\nu \right\rangle_{\Gamma(t)} = 0
\quad \forall\ \chi \in H^1(\Gamma(t))\,, \label{eq:ns_weakc} \\
& \left\langle \varkappa\,\vec\nu, \vec\eta \right\rangle_{\Gamma(t)}
+ \left\langle \nabs\,\vec \id, \nabs\,\vec \eta \right\rangle_{\Gamma(t)}
= 0 \quad \forall\ \vec\eta \in [H^1(\Gamma(t))]^d\,\label{eq:ns_weakd}
\end{align}
\end{subequations}
holds for almost all times $t \in (0,T]$. We notice that if $p \in \pspace$ is
part of a solution to (\ref{eq:ns_full_momentum}--g), then so is $p + c$ for an
arbitrary $c\in \R$.

An alternative to the weak formulation (\ref{eq:ns_weaka}--d) can be obtained
by
following the paper \cite{fluidfbp}. These authors show that it is possible
to derive an energy bound not only on the continuous level but also
for a discretization based on this formulation.
Due to the resulting structure of the equation,
we will refer to this form as the antisymmetric weak formulation,
and it is given as follows.
Given $\Gamma(0) = \Gamma_0$ and $\vec u(\cdot,0) = \vec u_0$, for
almost all $t\in(0,T)$ find $\Gamma(t)$ and $(\vec u, p, \varkappa)(t) \in
\uspace g \times \pnormspace \times H^1(\Gamma(t))$ such that
\begin{align}
& \tfrac{1}{2}\bigg[ \ddt (\rho\,\vec u, \vec \xi) + (\rho\,\vec u_t, \vec \xi)
+ (\rho, [(\vec u\,.\,\nabla)\,\vec u]\,.\,\vec \xi
- [(\vec u\,.\,\nabla)\,\vec \xi]\,.\,\vec u)\bigg] \nonumber \\
& \qquad +2\left(\mu\,\mat D(\vec u), \mat D(\vec \xi)\right)
- \left(p, \nabla\,.\,\vec \xi\right)
- \gamma\,\left\langle \varkappa\,\vec\nu, \vec\xi
\right\rangle_{\Gamma(t)}
= \left(\vec f, \vec \xi\right)\quad \forall\ \vec\xi \in \uspace 0 \,,
\label{eq:ns_weaka_antisym}
\end{align}
holds, together with (\ref{eq:ns_weakb}--d). For a derivation of
(\ref{eq:ns_weaka_antisym}) we refer to \cite{fluidfbp}, see also
\cite{Agnese}.

\subsection{ALE model and weak formulation} \label{sec:ale_model}
The Arbitrary Lagrangian Eulerian (ALE) technique consists in
reformulating the time derivative $\vec u_t$ by expressing it with respect to a
fixed reference configuration. A special homeomorphic map, called the ALE map,
associates, at each time $t$, a point in the current computational domain
$\Omega(t)$ to a point in the reference domain $\D$. Note that for convenience
we consider the domain $\Omega$ to be time-dependent, even though in our
situation it is fixed.

We now reformulate the two-phase Navier--Stokes system
(\ref{eq:ns_full_momentum}--g), which is expressed in terms of Eulerian
coordinates $\vec z\in\Omega(t)$, by rewriting the velocity time derivative
$\vec u_t$ with respect to the so called ALE coordinate $\vec q \in \D$.
Analogously to what is done in a front tracking approach to parametrize the
interface $\Gamma(t)$, we can define the fixed reference manifold $\D$ and
extend the map $\vec x(\cdot,t)$ to parametrize $\Omega(t)$ as
$\Omega(t)=\vec x(\D,t)$. This extended map clearly still satisfies $\Gamma(t)
= \vec x(\Upsilon,t)$, given that $\Upsilon\subset\D$. Moreover, we let
$\partial\D$ be partitioned as $\partial\D=\partial_1\D \cup \partial_2\D$
with $\partial_1\D \cap \partial_2\D = \emptyset$.

Now, let $h:\Omega(t)\times [0,T]\to \R$ be a function defined on the Eulerian
frame. The corresponding function on the ALE frame $\hat h$ is defined as
\begin{equation}\label{eq:h}
\hat h:\D\times [0,T]\to \R,\quad \hat h(\vec q,t)=h(\vec x(\vec q,t),t).
\end{equation}
In order to compute the time derivative of (\ref{eq:h}) with respect to the ALE
frame, using the chain rule, we have
$\frac{\partial\hat h(\vec q,t)}{\partial t}=\frac{\partial h(\vec x(\vec
q,t),t)}{\partial t}=\frac{\partial h(\vec z,t)}{\partial t}+\vec x_t(\vec q,t)
\,.\, \nabla\, h(\vec z,t)$. Therefore it holds that
$\frac{\partial h(\vec z,t)}{\partial t} =
\frac{\partial\hat h(\vec q,t)}{\partial t}-
\vec x_t(\vec q,t) \,.\, \nabla\, h(\vec z,t)$.
Finally, introducing the domain velocity
\begin{equation} \label{eq:W}
\W(\vec z, t) := \vec x_t(\vec q, t) \quad \forall\ \vec z = \vec x(\vec q,t)
\in \Omega(t)
\end{equation}
and the time derivative in the ALE frame
\begin{equation} \label{eq:ale_derivative}
\left.\frac{\partial h(\vec z,t)}{\partial t}\right|_{\D}:=
\frac{\partial\hat h(\vec q,t)}{\partial t} \quad
\forall\ \vec z = \vec x(\vec q,t) \in \Omega(t)\,,
\end{equation}
we obtain
\begin{equation}\label{eq:ale_dt}
h_t =\left.h_t\right|_{\D} -\vec{\mathcal{W}}\,.\,\nabla\, h\,.
\end{equation}
The identity (\ref{eq:ale_dt}) is naturally extended to vector valued functions.
We stress the fact that the domain velocity $\W$ for the interface points is
consistent with the interface velocity $\V$. Therefore it holds that
$\W\!\mid_{\Gamma(t)}=\V$.
We also point out that the ALE mapping is somehow arbitrary, apart from the
requirement of conforming to the evolution of the domain boundary. Indeed the
map of the boundary $\partial_i \D$ of the reference domain has then to
provide, at all $t$, the boundary $\partial_i \Omega(t)$ of the current
configuration, with $i=1,\,2$.

Using (\ref{eq:ale_dt}) in the momentum equation (\ref{eq:ns_full_momentum}),
we can rewrite it as:
\begin{align}
\rho\,(\left.\vec u_t\right|_{\D} +((\vec u - \W)\,.\, \nabla)\, \vec u)
 -2 \mu\,\nabla\,.\,\mat D(\vec u)+ \nabla\,p & = \vec f
\quad &&\mbox{in } \Omega_\pm(t)\,. \label{eq:ns_full_momentum_ale}
\end{align}
We can notice that, with respect to the original formulation, there is a
convective-type term due to the domain movement and the time derivative is
computed in the fixed reference frame $\D$. Moreover, contrary to the original
system, not only the regions $\Omega_\pm(t)$ are time dependent but also the
whole domain $\Omega(t)$ is time dependent.

Obviously, if the domain is fixed, the additional convective term is zero
and the time derivative in the ALE frame coincides with the usual time
derivative in the Eulerian frame. This means that $\W=\vec 0$ corresponds to a
pure Eulerian method, while $\W=\vec u$ corresponds to a fully Lagrangian
scheme. In our case, since we have a fixed domain $\Omega$, the ALE formulation
might not seem very useful. But, at the discrete level, we are concerned with
the evolution of the discrete triangulated domains. Therefore the discrete ALE
map describes the evolution of the grid during the domain movement. It is
indeed at the discrete level that the advantage of the ALE formulation emerges,
as in an ALE setting the time advancing scheme provides directly the evolution
of the unknowns at mesh nodes and thus that of the degrees of freedom of the
discrete solutions.

In order to define the ALE weak formulation, we need to use a different
functional setting for the test functions with respect to the one used for the
Eulerian weak formulations. The reason is that here it
needs to be defined on the moving domain $\Omega(t)$. Therefore we introduce
the admissible spaces of test functions on the reference domain $\D$
\begin{align*}
\uspaceref b &:= \{\vec \phi\in[H^1(\D)]^d:
\vec \phi =\vec b\quad \mbox{on }\partial_1\D\,,
\quad \vec \phi\,.\,\unitn=0 \quad \mbox{on }\partial_2\D\}\,,\\
\pspaceref &:= L^2(\D)\,,\qquad
\pnormspaceref := \{\eta\in\pspaceref : \int_\D\eta\dL{d}=0 \}\,,
\end{align*}
for a given $\vec b \in [H^1(\D)]^d$. Then, using the ALE mapping, we can
define the admissible spaces of test functions on the moving domain $\Omega(t)$
by setting
\begin{align*}
\uspaceale b &:= \{ \vec\phi:
\bigcup_{t \in [0,T]} \Omega(t)\times \{t\} \to \R^d,\quad
\vec \phi=\vec{\hat\phi}\circ\vec x^{\,-1},\quad\vec{\hat\phi}(\cdot,t)
\in\uspaceref b \quad \forall\ t \in [0,T]\}\,,\\
\pnormspaceale &:= \{ \eta:
\bigcup_{t \in [0,T]} \Omega(t)\times \{t\} \to \R,\quad
\eta=\hat\eta\circ\vec x^{\,-1},\quad\hat\eta(\cdot,t)\in\pnormspaceref
\quad \forall\ t \in [0,T] \}\,.
\end{align*}
Hence, the ALE weak formulation of (\ref{eq:ns_full_momentum_ale}),
(\ref{eq:ns_full_mass}--g)
is given as follows. Given $\Gamma(0) = \Gamma_0$ and
$\vec u(\cdot,0) = \vec u_0$,
find $(\Gamma(t))_{t\in[0,T]}$ and $(\vec u, p, \varkappa)\in \uspaceale g
\times \pnormspaceale \times H^1(\bigcup_{t\in[0,T]} \Gamma(t) \times \{t\})$
such that
\begin{subequations}
\begin{align}
& \left(\rho\,\left.\vec u_t\right|_{\D}, \vec \xi\right)_{\Omega(t)} +
\left(\rho\,((\vec u - \W)\,.\, \nabla )\,\vec u,\vec \xi\right)_{\Omega(t)}
+ 2\left(\mu\,\mat D(\vec u), \mat D(\vec \xi)\right)_{\Omega(t)} \nonumber \\
& \qquad - \left(p, \nabla\,.\,\vec \xi\right)_{\Omega(t)}
- \gamma\,\left\langle \varkappa\,\vec\nu, \vec\xi\right\rangle_{\Gamma(t)}
= \left(\vec f, \vec \xi\right)_{\Omega(t)}\quad \forall\ \vec\xi \in
\uspaceale 0 \,, \label{eq:ns_weaka_ale}\\
& \left(\nabla\,.\,\vec u, \varphi\right)_{\Omega(t)} = 0
\quad \forall\ \varphi \in \pnormspaceale\,, \label{eq:ns_weakb_ale}
\end{align}
\end{subequations}
and (\ref{eq:ns_weakc},d) hold for almost all times $t \in (0,T]$.
Here $(\cdot,\cdot)_{\Omega(t)}$ denotes the $L^2$--inner products on
$\Omega(t)$.

\setcounter{equation}{0}
\section{Finite element approximations}\label{sec:ns_fem}
\subsection{Eulerian finite element approximations}\label{sec:ns_fem_antisym}
We consider the partitioning $t_m =m\,\tau$, $m=0,\ldots, M$, of $[0,T]$
into uniform time
steps $\tau=\frac{T}{M}$, see \cite{fluidfbp}.
Moreover, let ${\cal T}^m$, $\forall\ m\ge 0$, be a regular partitioning of the
domain $\Omega$ into disjoint open simplices $\sigmaO^m_j$, $j = 1 ,\ldots,
J^m_\Omega$. From now on, the domain $\Omega$ which we consider is the
polyhedral domain defined by the triangulation ${\cal T}^m$. On ${\cal T}^m$ we
define the finite element spaces
\begin{equation*}
S^m_k := \{\chi \in C(\overline{\Omega}) : \chi\!\mid_{\sigmaO^m}
\in \mathcal{P}_k(\sigmaO^m) \ \forall\ \sigmaO^m \in {\cal T}^m\}\,,
\quad k \in \mathbb{N}\,,
\end{equation*}
where $\mathcal{P}_k(\sigmaO^m)$ denotes the space of polynomials of degree $k$
on $\sigmaO^m$. Moreover, $S^m_0$ is the space of piecewise constant functions
on ${\cal T}^m$ and let $\vec I^m_k$ be the standard interpolation operator
onto $[S^m_k]^d$.

Let $\uspacedisc{g}{m}\subset\uspacesimple(\vec I_k^m\vec g)$ and
$\pspace^m\subset\pspace$ be the finite element spaces we use for the
approximation of velocity and pressure, and let $\pnormspace^m:= \pspace^m \cap
\pnormspace$. In this paper, we restrict ourselves to the choice
P2--(P1+P0), i.e.\ we set
$\uspacedisc{0}{m}=[S^m_2]^d\cap\uspace{0}$ and
$\pspace^m = S^m_1+S^m_0$, but generalizations to other spaces are easily
possible.

We consider a fitted finite element approximation for the evolution of the
interface $\Gamma(t)$. Let $\Gamma^m\subset\R^d$ be a $(d-1)$-dimensional
polyhedral surface approximating the closed surface $\Gamma(t_m)$, $m=0
,\ldots, M$. Let $\Omega^m_+$ denote the exterior of $\Gamma^m$ and let
$\Omega^m_-$ be the interior of $\Gamma^m$, where we assume that $\Gamma^m$ has
no self-intersections. Then $\Omega = \Omega_-^m \cup \Gamma^m \cup
\Omega_+^m$, and the fitted nature of
our method implies that
\begin{equation*}
\overline{\Omega^m_+} = \bigcup_{o \in \mathcal{T}^m_+} \overline{o}
\quad\text{and}\quad
\overline{\Omega^m_-} = \bigcup_{o \in \mathcal{T}^m_-} \overline{o} \,,
\end{equation*}
where $\mathcal{T}^m = \mathcal{T}^m_+ \cup \mathcal{T}^m_-$ and
$\mathcal{T}^m_+ \cap \mathcal{T}^m_- = \emptyset$. Let $\vec \nu^m$ denote the
piecewise constant unit normal to $\Gamma^m$ such that $\vec\nu^m$ points into
$\Omega^m_+$.
In addition, let
$\mu^m = \mu_+\,\charfcn{\Omega^m_+} + \mu_-\,\charfcn{\Omega^m_-}\in S^m_0$
and
$\rho^m = \rho_+\,\charfcn{\Omega^m_+} + \rho_-\,\charfcn{\Omega^m_-}
\in S^m_0$.

In order to define the parametric finite element spaces on $\Gamma^m$, we
assume that $\Gamma^m=\bigcup_{j=1}^{J_\Gamma} \overline{\sigma^m_j}$, where
$\{\sigma^m_j\}_{j=1}^{J_\Gamma}$ is a family of mutually disjoint open
$(d-1)$-simplices with vertices $\{\vec q^m_k\}_{k=1}^{K_\Gamma}$. Then
we define $\Vh := \{\vec\chi \in [C(\Gamma^m)]^d:\vec\chi\!\mid_{\sigma^m_j}
\in \mathcal{P}_1(\sigma^m_j), j=1,\ldots, J_\Gamma\} =: [\Wh]^d$, where $\Wh
\subset H^1(\Gamma^m)$ is the space of scalar continuous piecewise linear
functions on $\Gamma^m$, with $\{\chi^m_k\}_{k=1}^{K_\Gamma}$ denoting the
standard basis of $\Wh$. As usual, we parametrize the new surface
$\Gamma^{m+1}$ over $\Gamma^m$ using a parametrization $\vec X^{m+1} \in \Vh$,
so that $\Gamma^{m+1} = \vec X^{m+1}(\Gamma^m)$. Moreover, let
$\langle\cdot,\cdot\rangle_{\Gamma^m}^h$ be the mass lumped inner product on
$\Gamma^m$ defined as
\begin{equation} \label{eq:masslump}
\left\langle v, w \right\rangle^h_{\Gamma^m} :=
\tfrac1d \sum_{j=1}^{J_\Gamma} \surfvol(\sigma^m_j)\,
\sum_{k=1}^d (v\,w)((\vec q^m_{j_k})^-),
\end{equation}
where $\{\vec q^m_{j_k}\}_{k=1}^d$ are the vertices of $\sigma^m_j$, and
where we define the limit $v((\vec q^m_{j_k})^-)
:= \underset{\sigma^m_j\ni \vec p\to \vec q^m_{j_k}}{\lim}\, v(\vec p)$. We
naturally extend (\ref{eq:masslump}) to vector valued functions. Finally, let
$\langle\cdot,\cdot\rangle_{\Gamma^m}$ denote the standard $L^2$--inner product
on $\Gamma^m$.

Our fully discrete finite element approximations based on the standard weak
formulation (\ref{eq:ns_weaka}--d) are given as follows.
Here we consider an explicit and an implicit treatment of the convective term
in the Navier--Stokes equation.

$(\schemeAex)$: Let
$\Gamma^0$, an approximation to $\Gamma(0)$, and $\vec U^0\in \uspacedisc{g}{0}$
be given. For $m=0,\ldots, M-1$, find $(\vec U^{m+1},P^{m+1}, \vec X^{m+1},
\kappa^{m+1}) \in \uspacedisc{g}{m}\times \pnormspace^m \times \Vh \times \Wh$
such that
\begin{subequations}
\begin{align}
& \left( \rho^m\,\frac{\vec U^{m+1} - \vec I^m_2\,\vec U^m}{\tau}, \vec
\xi \right) + \left(\rho^m(\vec I^m_2\vec U^m\,.\,\nabla)\,\vec U^{m+1}\,,
\,\vec \xi \right)
+ 2\left(\mu^m\,\mat D(\vec U^{m+1}), \mat D(\vec \xi) \right)
\nonumber \\ & \quad
- \left(P^{m+1}, \nabla\,.\,\vec \xi\right)
- \gamma\,\left\langle \kappa^{m+1}\,\vec\nu^m ,
\vec\xi\right\rangle_{\Gamma^m}
= \left(\rho^m\,\vec f_1^{m+1}+\vec f_2^{m+1}, \vec \xi\right)
\quad \forall\ \vec\xi \in \uspacedisc{0}{m}\,, \label{eq:ns_HGa} \\
& \left(\nabla\,.\,\vec U^{m+1}, \varphi\right) = 0
\quad \forall\ \varphi \in \pnormspace^m\,,\label{eq:ns_HGb} \\
& \left\langle \frac{\vec X^{m+1} - \vec\id}{\tau} ,\chi\,\vec\nu^m
\right\rangle_{\Gamma^m}^h - \left\langle \vec U^{m+1}, \chi\,\vec\nu^m
\right\rangle_{\Gamma^m} = 0 \quad \forall\ \chi \in \Wh\,, \label{eq:ns_HGc}\\
& \left\langle \kappa^{m+1}\,\vec\nu^m, \vec\eta \right\rangle_{\Gamma^m}^h
+ \left\langle \nabs\,\vec X^{m+1}, \nabs\,\vec \eta \right\rangle_{\Gamma^m} =
0 \quad \forall\ \vec\eta \in \Vh \label{eq:ns_HGd}
\end{align}
\end{subequations}
and set $\Gamma^{m+1} = \vec X^{m+1}(\Gamma^m)$. Here we have defined
$\vec f_i^{m+1}(\cdot) := \vec I^m_2\,\vec f_i(\cdot,t_{m+1})$, $i=1,2$.

$(\schemeAim)$: Let
$\Gamma^0$, an approximation to $\Gamma(0)$, and $\vec U^0\in \uspacedisc{g}{0}$
be given. For $m=0,\ldots, M-1$, find $(\vec U^{m+1},P^{m+1}, \vec X^{m+1},
\kappa^{m+1}) \in \uspacedisc{g}{m}\times \pnormspace^m \times \Vh \times \Wh$
such that
\begin{align}
& \left( \rho^m\,\frac{\vec U^{m+1} - \vec I^m_2\,\vec U^m}{\tau}, \vec
\xi \right) + \left(\rho^m( \vec U^{m+1}\,.\,\nabla)\,\vec U^{m+1}\,,
\,\vec \xi \right)
+ 2\left(\mu^m\,\mat D(\vec U^{m+1}), \mat D(\vec \xi) \right)
\nonumber \\ & \qquad
- \left(P^{m+1}, \nabla\,.\,\vec \xi\right)
- \gamma\,\left\langle \kappa^{m+1}\,\vec\nu^m ,
\vec\xi\right\rangle_{\Gamma^m}
= \left(\rho^m\,\vec f_1^{m+1}+\vec f_2^{m+1}, \vec \xi\right)
\quad \forall\ \vec\xi \in \uspacedisc{0}{m}\,, \label{eq:ns_HGimpa}
\end{align}
and (\ref{eq:ns_HGb}--d), and set $\Gamma^{m+1} = \vec X^{m+1}(\Gamma^m)$.
Since the convective term in (\ref{eq:ns_HGa}) is treated explicitly,
the scheme $(\schemeAex)$ is a linear scheme that leads to a coupled linear
system of equations for the unknowns
$(\vec U^{m+1}, P^{m+1}, \vec X^{m+1}, \kappa^{m+1})$ at each time level.
On the other hand, in (\ref{eq:ns_HGimpa}) the convective term is treated
implicitly, and so $(\schemeAim)$ is a nonlinear scheme that leads to a coupled
nonlinear system of equations for the unknowns $(\vec U^{m+1}, P^{m+1}, \vec
X^{m+1}, \kappa^{m+1})$ at each time level.

An alternative fully discrete finite element approximation is
based on the formulation (\ref{eq:ns_weaka_antisym}), (\ref{eq:ns_weakb}--d).
In particular, we consider the following scheme.

$(\schemeB)$: Let $\Gamma^0$, an approximation to $\Gamma(0)$, and $\vec U^0\in
\uspacedisc{g}{0}$ be given. For $m=0,\ldots, M-1$, find $(\vec U^{m+1},P^{m+1},
\vec X^{m+1}, \kappa^{m+1}) \in \uspacedisc{g}{m}\times \pnormspace^m \times
\Vh \times \Wh$ such that
\begin{align}
& \frac1\tau \left(\tfrac{1}{2}(I^m_0\,\rho^{m-1}+\rho^m)\,\vec U^{m+1} -
I^m_0\,\rho^{m-1} \vec I^m_2\,\vec U^m, \vec \xi \right)
+ \tfrac{1}{2}\left((\rho^m\,\vec I^m_2\,\vec U^m\,.\,\nabla)\,
\vec U^{m+1}\,,\,\vec \xi \right)
\nonumber \\ & \qquad
- \tfrac{1}{2} \left((\rho^m\,
\vec I^m_2 \, \vec U^m\,.\,\nabla)\,\vec \xi\,,\,\vec U^{m+1} \right)
+ 2\left(\mu^m\,\mat D(\vec U^{m+1}), \mat D(\vec \xi) \right)
- \left(P^{m+1}, \nabla\,.\,\vec \xi\right) \nonumber \\
& \qquad - \gamma\,\left\langle
\kappa^{m+1}\,\vec\nu^m,\vec\xi\right\rangle_{\Gamma^m}
= \left(\rho^m\,\vec f_1^{m+1}+\vec f_2^{m+1}, \vec \xi\right)
\quad \forall\ \vec\xi \in \uspacedisc{0}{m}\,,\label{eq:ns_HGa_antisym}
\end{align}
and (\ref{eq:ns_HGb}--d), and set $\Gamma^{m+1} = \vec X^{m+1}(\Gamma^m)$.
We note that the scheme $(\schemeB)$ corresponds to
\cite[(4.6a--d)]{fluidfbp} in the setting of an unfitted finite element
approximation. Similarly to \cite[Theorem~4.1]{fluidfbp} it is then
straightforward to show that there exists a unique solution to
$(\schemeB)$ that fulfils an energy inequality. In the unfitted case this
energy inequality leads to an overall stability result if the bulk mesh does
not change. Of course, in the fitted case considered here this will never be
the case, unless the interface is stationary. Nevertheless, it is of interest
to see how the scheme $(\schemeB)$ performs in practice, compared to
$(\schemeAex)$ and $(\schemeAim)$.

Moreover, we remark that all three schemes, $(\schemeAex)$, $(\schemeAim)$
and $(\schemeB)$, share the equations (\ref{eq:ns_HGd}). Apart from
defining the discrete curvature, the equations (\ref{eq:ns_HGd}) can also
be viewed as a side constraint on the distribution of mesh points,
see \cite{triplej,gflows3d} for further details. In particular, they lead
to an implicit tangential motion of the vertices on the discrete interface,
so that the vertices of
$\Gamma^m$ will in general be very well distributed in practice.
For example, in the case $d=2$ it is possible to prove an
asymptotic equidistribution property, see \cite[Theorem~5]{Agnese}.

Finally, we stress that the three schemes $(\schemeAex)$, $(\schemeAim)$
and $(\schemeB)$, all feature at least one term containing the
interpolated velocity $\vec I^m_2\,\vec U^m$, and so an interpolation
of the old velocity $\vec U^m$ onto the new mesh $\mathcal{T}^m$ is necessary,
see also \S\ref{sec:velocity_interpolation} below.
In order to avoid such an interpolation, we consider an ALE formulation next.

\subsection{ALE finite element approximation}\label{sec:ale_fem}
In order to derive our fully discrete ALE finite element approximation, we
first introduce a semidiscrete continuous-in-time fitted finite element
approximation of (\ref{eq:ns_weaka_ale}--d).
Let ${\mathfrak T}^h$ be a regular
partitioning of the reference domain $\D$ into disjoint open simplices
$\sigmaO^h_j$, $j = 1 ,\ldots, J^h_\D$. From now on, the reference domain $\D$,
which we consider, is the polyhedral domain defined by the triangulation
${\mathfrak T}^h$. On ${\mathfrak T}^h$ we define the finite element spaces
\begin{equation*}
S^h_k(\D) := \{\chi \in C(\overline{\D}) : \chi\!\mid_{\sigmaO^h}
\in \mathcal{P}_k(\sigmaO^h) \ \forall\ \sigmaO^h \in {\mathfrak T}^h\}\,,
\quad k \in \mathbb{N}\,,
\end{equation*}
as well as $S^h_0(\D)$, the space of piecewise constant
functions on ${\mathfrak T}^h$.

Then, using the discrete ALE mapping
$\vec x_h$, we can define the finite element spaces on
the discretized moving domains $\Omega^h(t) = \vec x_h(\D,t)$ by
\begin{equation*}
S^h_k(\Omega^h(t)) :=
\{ \eta: \Omega^h(t)\to \R,\quad
\eta=\hat\eta\circ\vec x_h^{\,-1},\quad\hat\eta\in S^h_k(\D) \}\,,\quad
k \in \mathbb{N} \cup \{0\}\,.
\end{equation*}
Here $\vec x_h\in [S^h_1]^d$ maps each simplex with straight faces in
${\mathfrak T}^h$ to a simplex with straight faces belonging to the
triangulation $\mathcal{T}^h(t)$ of the discretized moving domain
$\Omega^h(t)$.
For later use, we let
$\uspacesemidiscale{g}{h}{t} \subset [S^h_2(\Omega^h(t))]^d$ and
$\pnormspaceale^h(t) \subset S^h_1(\Omega^h(t)) + S^h_0(\Omega^h(t))$
denote the appropriate spaces for the discrete velocities and pressure.
In addition, we also define the discrete domain velocity as
\begin{equation}
\W^h(\vec z, t):=
\sum_{k=1}^{K_\Omega}\left[\ddt\,\vec p^h_k(t)\right] \varphi^h_k(\vec z, t)
\in [S^h_1(\Omega^h(t))]^d\,,
\end{equation}
with $\{\varphi^h_k(\cdot,t)\}_{k=1}^{K_\Omega}$ denoting the standard basis of
$S^h_1(\Omega^h(t))$,
and $\{\vec p^h_k(t)\}_{k=1}^{K_\Omega}$ denoting the vertices of
$\mathcal {T}^h(t)$.

Given $(\Gamma^h(t))_{t\in[0,T]}$, a family of
$(d-1)$-dimensional polyhedral surfaces with unit normal $\vec\nu^h(t)$
pointing into $\Omega^h_+(t)$, the exterior of $\Gamma^h(t)$,
analogously to the fully discrete case discussed in the previous section, we
introduce the
piecewise linear finite element spaces $\Wht$ and $\Vht$, with
$\{\chi^h_k(\cdot,t)\}_{k=1}^{K_\Gamma}$ denoting the standard basis of the
former. Hence $\chi^h_k(\vec q^h_l(t),t) = \delta_{kl}$ for all
$k,l \in \{1,\ldots,K_\Gamma\}$ and $t \in [0,T]$, where
$\{\vec q^h_k(t)\}_{k=1}^{K_\Gamma}$ are the vertices of $\Gamma^h(t)$. We also
notice that the discrete interface velocity $\V^h$ defined as
\begin{equation}\label{eq:Vsemidisc}
\V^h(\vec z, t):=
\sum_{k=1}^{K_\Gamma}\left[\ddt\,\vec q^h_k(t)\right] \chi^h_k(\vec z, t)
\in \Vht\,,
\end{equation}
see e.g. \cite[(3.3)]{tpfs}, is simply the trace of $\W^h$ on $\Gamma^h$.
We also introduce the interior $\Omega^h_-(t) = \Omega(t) \setminus
\overline{\Omega^h_+(t)}$, and the $L^2$--inner products
$(\cdot,\cdot)_{\Omega^h(t)}$ and $\langle\cdot,\cdot\rangle_{\Gamma^h(t)}$
on $\Omega^h(t)$ and $\Gamma^h(t)$, respectively. Moreover, let
$\langle\cdot,\cdot\rangle_{\Gamma^h(t)}^h$ be the mass lumped inner product on
$\Gamma^h(t)$.
Finally, we let $\vec f_i^h(\cdot,t) := \vec I^h_2\,\vec f_i(\cdot,t)$,
$i=1,2$, where $\vec I^h_2$ denotes the standard interpolation operator onto
$[S^h_2(\Omega^h(t))]^d$, and set
$\mu^h(t) = \mu_+\,\charfcn{\Omega^h_+(t)} + \mu_-\,\charfcn{\Omega^h_-(t)}\in
S^h_0(\Omega^h(t))$
and
$\rho^h(t)
= \rho_+\,\charfcn{\Omega^h_+(t)} + \rho_-\,\charfcn{\Omega^h_-(t)}\in
S^h_0(\Omega^h(t))$.

Then, given $\Gamma^h(0)$ and $\vec U^h(0)\in \uspacesemidiscale{g}{h}{0}$, for
$t\in (0,T]$ find $\Gamma^h(t)$ and $(\vec U^h, P^h,
\V^h, \kappa^h)(t) \in \uspacesemidiscale{g}{h}{t} \times \pnormspaceale^h(t)
\times \Vht \times \Wht$ such that
\begin{subequations}
\begin{align}
& \left(\rho^h\,\left.\vec U^h_t\right|_{\D}, \vec \xi \right)_{\Omega^h(t)} +
\left(\rho^h((\vec U^h-\W^h)\,.\,\nabla)\,\vec U^h\,,
\,\vec \xi \right)_{\Omega^h(t)}
+2\left(\mu^h\,\mat D(\vec U^h), \mat D(\vec \xi) \right)_{\Omega^h(t)}
\nonumber \\ & \qquad
- \left(P^h, \nabla\,.\,\vec \xi\right)_{\Omega^h(t)}
- \gamma\,\left\langle \kappa^h\,\vec\nu^h,
\vec\xi\right\rangle_{\Gamma^h(t)} =
\left(\rho^h\,\vec f_1^h+\vec f_2^h, \vec \xi\right)_{\Omega^h(t)}
\quad \forall\ \vec\xi \in \uspacesemidiscale{0}{h}{t}\,,
\label{eq:ns_semi_alea}\\
& \left(\nabla\,.\,\vec U^h, \varphi\right)_{\Omega^h(t)} = 0
\quad \forall\ \varphi \in \pnormspaceale^h(t)\,, \label{eq:ns_semi_aleb} \\
& \left\langle \V^h , \chi\,\vec\nu^h
\right\rangle_{\Gamma^h(t)}^h - \left\langle \vec U^h, \chi\,\vec\nu^h
\right\rangle_{\Gamma^h(t)} = 0 \quad \forall\ \chi \in \Wht\,,
\label{eq:ns_semi_alec} \\
& \left\langle \kappa^h\,\vec\nu^h, \vec\eta \right\rangle_{\Gamma^h(t)}^h
+ \left\langle \nabs\,\vec\id, \nabs\,\vec \eta \right\rangle_{\Gamma^h(t)} = 0
\quad \forall\ \vec\eta \in \Vht\,. \label{eq:ns_semi_aled}
\end{align}
\end{subequations}
The semidiscrete continuous-in-time fitted finite element scheme
(\ref{eq:ns_semi_alea}--d) is a system of ODEs. Indeed, on defining the
appropriate time-dependent matrices and vectors, we can rewrite
(\ref{eq:ns_semi_alea}--d) in algebraic form as
\begin{subequations}
\begin{align}
& \vec M^h_\Omega(t)\ddt \vec U + \vec B^h_\Omega(t) \vec U +
\vec D^h_\Omega\big(t,\vec U^h,\W^h\big)\vec U
+ \vec C^h_\Omega(t) P
- \gamma \Nbulk^h(t) \kappa = \vec c^h(t,\vec f_1^h,\vec f_2^h)\,,
\label{eq:ns_ode_alea} \\
& [\vec C^h_\Omega(t) ]^T \vec U = 0\,, \label{eq:ns_ode_aleb} \\
& [\vec N_\Gamma^h(t)]^T \ddt \vec X - \Nbulk^h(t)\vec U = 0\,,
\label{eq:ns_ode_alec} \\
& \vec N_\Gamma^h(t) \kappa + \vec A^h_\Gamma(t)\vec X = 0\,,
\label{eq:ns_ode_aled}
\end{align}
\end{subequations}
where $\vec U$, $P$, $\kappa$ and $\vec X$ are the unknown vectors for the
velocity, pressure, interface curvature and interface position, respectively.
See \cite[(6.15)]{Agnese} for details.

It is now straightforward to obtain a fully discrete finite element
approximation of (\ref{eq:ns_ode_alea}--d), since the only missing part is the
temporal discretization. On recalling the partitioning $t_m =m\,\tau$,
$m=0,\ldots, M$, of $[0,T]$, and on applying a semi-implicit Euler scheme to
(\ref{eq:ns_ode_alea}--d), we obtain
\begin{subequations}
\begin{align}
& \frac{1}{\tau}\vec M^h_\Omega(t^{m+1})\vec U^{m+1} +
\vec B^h_\Omega(t^m) \vec U^{m+1} +
\vec D^h_\Omega\big(t^m,\vec U^{m+1},\W^{m+1}\big)\vec U^{m+1}
+ \vec C^h_\Omega(t^m) P^{m+1}
\nonumber \\ & \qquad
- \gamma \Nbulk^h(t^m) \kappa^{m+1}
= \frac{1}{\tau} \vec M^h_\Omega(t^{m+1})\vec U^m
+ \vec c^h(t^m,\vec f_1^{m+1},\vec f_2^{m+1})\,,\label{eq:ns_ode_disc_alea} \\
& [\vec C^h_\Omega(t^m) ]^T \vec U^{m+1} = 0\,,
\label{eq:ns_ode_disc_aleb} \\
& \frac{1}{\tau}[\vec N_\Gamma^h(t^m)]^T \vec X^{m+1}
- \Nbulk^h(t^m)\vec U^{m+1} = \frac{1}{\tau}[\vec N_\Gamma^h(t^m)]^T
\vec X^m\,, \label{eq:ns_ode_disc_alec}\\
& \vec N_\Gamma^h(t^m) \kappa^{m+1} + \vec A^h_\Gamma(t^m)\vec X^{m+1}
 = 0\,, \label{eq:ns_ode_disc_aled}
\end{align}
\end{subequations}
where $\W^m$, $m=1,\ldots,M$ are some fully discrete bulk mesh velocities
which, starting from $\Omega^0$, induce a sequence of bulk meshes
$\Omega^m=\vec x^m(\D)$, $m=1,\ldots,M$, where
$\vec x^m \in [S^h_1(\Omega^m)]^d$
is an approximation to $\vec x^h(\cdot,t_m)$.

Following \cite{NobilePhd}, we now rewrite (\ref{eq:ns_ode_disc_alea}--d)
into a form that is common in the finite element literature.
Firstly, let $\uspacediscale{g}{m}= \{ \vec\phi:
\Omega^m \to \R^d,\
\vec \phi=\vec{\hat\phi}\circ(\vec x^m)^{\,-1},\ \vec{\hat\phi}
\in \mathcal{U}(\vec g) \cap [S^h_2(\D)]^d\}$
and $\pnormspaceale^m=
\{ \eta: \Omega^m\to \R,\
\eta=\hat\eta\circ(\vec x^m)^{\,-1},\ \hat\eta\in\pnormspaceref
\cap (S^h_1(\D) + S^h_0(\D)) \}$.
Now to ease the notation, we introduce the following convention. For a function
$h^k$ defined on $\Omega^k$, we can immediately transport it to any other
configuration $\Omega^m$, with $k\neq m$, using the mapping
$\vec x^k\circ(\vec x^m)^{-1}$. Hence, when we need to integrate
$h^k$ on $\Omega^m$ we will write simply $\int_{\Omega^m} h^k \dL d$ instead of
$\int_{\Omega^m} h^k\circ \vec x^k\circ(\vec x^m)^{-1}\dL d$, and we will
adopt the same notation used for the inner products.
Then, our ALE finite element approximation (\ref{eq:ns_ode_disc_alea}--d),
which is based on the variational formulation
(\ref{eq:ns_weaka_ale}--d), can be written as follows.

$(\schemeALE)$: Let $\Omega^0$, $\Gamma^0$, and
$\vec U^0\in \uspacediscale{g}{0}$ be given.
For $m=0,\ldots, M-1$, find
\linebreak
$(\vec U^{m+1},P^{m+1}, \vec X^{m+1}, \kappa^{m+1},
\W^{m+1})
\in \uspacediscale{g}{m+1}\times \pnormspaceale^{m+1} \times \Vh \times \Wh
\times [S^h_1(\Omega^m)]^d$ such that
$\W^{m+1}\!\mid_{\Gamma^m} = \frac1\tau(\vec X^{m+1} -
\vec\id\!\mid_{\Gamma^m})$ and such that
\begin{subequations}
\begin{align}
& \left( \frac{\rho^m}{\tau}\,\vec U^{m+1}, \vec \xi \right)_{\Omega^{m+1}}
+ \left(\rho^m( (\vec U^{m+1}- \W^{m+1})\,.\,\nabla)\,\vec U^{m+1}\,,
\,\vec \xi \right)_{\Omega^m}\nonumber \\
& \qquad + 2\left(\mu^m\,\mat D(\vec U^{m+1}), \mat D(\vec \xi)
\right)_{\Omega^m} - \left(P^{m+1}, \nabla\,.\,\vec \xi\right)_{\Omega^m}
- \gamma\,\left\langle \kappa^{m+1}\,\vec\nu^m , \vec\xi\right\rangle_{\Gamma^m}
\nonumber \\
& \qquad= \left(\frac{\rho^m}{\tau}\,\vec U^m,\vec \xi \right)_{\Omega^{m+1}}
+ \left(\rho^m\,\vec f_1^{m+1}+\vec f_2^{m+1}, \vec \xi\right)_{\Omega^m}
\quad \forall\ \vec\xi \in \uspacediscale{0}{m}\,, \label{eq:ale_HGa} \\
& \left(\nabla\,.\,\vec U^{m+1}, \varphi\right)_{\Omega^m} = 0
\quad \forall\ \varphi \in \pnormspaceale^m\,,\label{eq:ale_HGb}
\end{align}
\end{subequations}
and (\ref{eq:ns_HGc},d), and set $\Gamma^{m+1} = \vec X^{m+1}(\Gamma^m)$,
as well as $\Omega^{m+1} = \vec x^{m+1}(\D)$, where
$\vec x^{m+1} = \vec x^m + \tau\,\W^{m+1}$. Here we note that
in practice the bulk mesh velocity $\W^{m+1}$ will be obtained with the
help of an iterative solution and smoothing procedure, see the next section for
more details.
Moreover, we recall that the side constraint (\ref{eq:ns_HGd}) will lead
to a uniform distribution of the vertices on $\Gamma^m$
in the case $d=2$, and to well distributed vertices in $d=3$.

We stress that since the
ALE formulation describes the evolution of the solution along trajectories, the
velocity $\vec U^{m}$ is defined on all the domains, and so in
particular on $\Omega^{m}$ and $\Omega^{m+1}$. Hence no interpolation
of the velocity is necessary in the ALE formulation, unless a full bulk
remeshing is performed because of severe mesh deformations.
Finally we remark that at the discrete level, the domain $\Omega^m$ often
plays the role of the reference manifold, so that $\D$ is never used
in practice.

\setcounter{equation}{0}
\section{Solution methods}\label{sec:solution_method}

\subsection{Eulerian schemes}

In this subsection we discuss the solution methods and other implementation
issues for the three Eulerian schemes $(\schemeAex)$, $(\schemeAim)$ and
$(\schemeB)$.

\subsubsection{Algebraic formulations} \label{sec:algebraic_system}
As is standard practice for the solution of linear systems arising from
discretizations of Navier--Stokes equations, we avoid the complications of the
constrained pressure space $\pnormspace^m$ in practice by considering an
overdetermined linear system with $\pspace^m$ instead.
For example, for the scheme $(\schemeAex)$ we consider the system
(\ref{eq:ns_HGa}--d) with $\pnormspace^m$ replaced by $\pspace^m$.
For non-zero Dirichlet boundary data $\vec g$ this requires a correction term
on the right hand side of (\ref{eq:ns_HGb}). In particular, we replace
(\ref{eq:ns_HGb}) with
\begin{equation} \label{eq:LAb}
 \left(\nabla\,.\,\vec U^{m+1}, \varphi\right) =
 \frac{\left(\varphi, 1\right)}{\vol(\Omega)}\, \int_{\partial_1\Omega}
(\vec I^m_2\,\vec g) \,.\, \unitn \dH{d-1} \quad \forall\ \varphi \in
\pspace^m\,.
\end{equation}
The linear system of equations arising at each time step of the scheme
$(\schemeAex)$ then can be formulated as: Find $(\vec U^{m+1},P^{m+1},
\kappa^{m+1},\delta\vec X^{m+1})$, where
$\vec X^{m+1} = \vec X^m+ \delta\vec X^{m+1}$, such that
\begin{equation}
\begin{pmatrix}
\vec B_\Omega & \vec C_\Omega & -\gamma\,\Nbulk & 0 \\
\vec C^T_\Omega & 0 & 0 & 0 \\
\NbulkT & 0 & 0 & -\frac1{\tau}\,\vec N_\Gamma^T \\
0 & 0 & \vec N_\Gamma & \vec A_\Gamma
\end{pmatrix}
\begin{pmatrix}
\vec U^{m+1} \\
P^{m+1} \\
\kappa^{m+1} \\
\delta\vec X^{m+1}
\end{pmatrix}
=
\begin{pmatrix}
\vec c \\
\vec \beta \\
0 \\
-\vec A_\Gamma\,\vec X^m
\end{pmatrix} \,,
\label{eq:ns_algebraic}
\end{equation}
where $(\vec U^{m+1},P^{m+1},\kappa^{m+1},\delta\vec X^{m+1})$ here,
in a slight abuse of notation, denote the
coefficients of these finite element functions with respect to some chosen
basis functions.
Moreover, $\vec X^m$ denotes the coefficients of
$\vec\id\!\mid_{\Gamma^m}$ with respect to the basis of $\Vh$. The definitions
of the matrices and vectors in (\ref{eq:ns_algebraic}) directly follow from
(\ref{eq:ns_HGa},c,d), (\ref{eq:LAb}),
and their exact definitions can be found in e.g.\ \cite[\S5.5]{Agnese}.
Note that
for the submatrices we have used the convention that the subscripts refer to
the test and trial domains, respectively. A single subscript is used where the
two domains are the same.
Of course, the algebraic formulation of the linear scheme
$(\schemeB)$ can also be stated in the form (\ref{eq:ns_algebraic}).

In order to solve the nonlinear scheme $(\schemeAim)$,
we employ the following fixed point iteration.
Let $\Gamma^m$ and $\vec U^m\in \uspacedisc{g}{m-1}$ be given.
Let $\vec U^{m+1,0}=\vec I^m_2\,\vec U^m$ and fix $\epsilon_f > 0$.
Then, for $s \geq 0$,
find $(\vec U^{m+1,s+1},P^{m+1,s+1}, \vec X^{m+1,s+1}, \kappa^{m+1,s+1}) \in
\uspacedisc{g}{m}\times \pspace^m \times \Vh \times \Wh$ such that
\begin{subequations}
\begin{align}
& \left( \rho^m\,\frac{\vec U^{m+1,s+1} - \vec I^m_2\,\vec U^m}{\tau}, \vec
\xi \right) + \left((\rho^m\,\vec U^{m+1,s}\,.\,\nabla)\,\vec U^{m+1,s+1}\,,
\,\vec \xi \right)
+ 2\left(\mu^m\,\mat D(\vec U^{m+1,s+1}), \mat D(\vec \xi) \right)
\nonumber \\ &\quad
- \left(P^{m+1,s+1}, \nabla\,.\,\vec \xi\right)
- \gamma\,\left\langle \kappa^{m+1,s+1}\,\vec\nu^m,\vec\xi\right\rangle_{\Gamma^m}
= \left(\rho^m\,\vec f_1^{m+1}+\vec f_2^{m+1}, \vec \xi\right)
\quad \forall\ \vec\xi \in \uspacedisc{0}{m}\,, \label{eq:ns_HGfixedpointa} \\
& \left(\nabla\,.\,\vec U^{m+1,s+1}, \varphi\right) =
 \frac{\left(\varphi, 1\right)}{\vol(\Omega)}\, \int_{\partial_1\Omega}
(\vec I^m_2\,\vec g) \,.\, \unitn \dH{d-1}
\quad \forall\ \varphi \in \pspace^m\,,\label{eq:ns_HGfixedpointb} \\
& \left\langle \frac{\vec X^{m+1,s+1} - \vec\id}{\tau} ,\chi\,\vec\nu^m
\right\rangle_{\Gamma^m}^h - \left\langle \vec U^{m+1,s+1}, \chi\,\vec\nu^m
\right\rangle_{\Gamma^m} = 0 \quad \forall\ \chi \in \Wh\,,
\label{eq:ns_HGfixedpointc}\\
& \left\langle \kappa^{m+1,s+1}\,\vec\nu^m, \vec\eta \right\rangle_{\Gamma^m}^h
+ \left\langle \nabs\,\vec X^{m+1,s+1}, \nabs\,\vec \eta
\right\rangle_{\Gamma^m} = 0 \quad \forall\ \vec\eta \in \Vh\,,
\label{eq:ns_HGfixedpointd}
\end{align}
\end{subequations}
and repeat the iteration until $\|\vec U^{m+1,s+1}-\vec U^{m+1,s}\|_{L^\infty}
\leq\epsilon_f$. Finally, set $\Gamma^{m+1}=$ \linebreak
$\vec X^{m+1,s+1}(\Gamma^m)$
and $(\vec U^{m+1}, P^{m+1}, \kappa^{m+1}) = (\vec U^{m+1,s+1}, P^{m+1,s+1},
\kappa^{m+1,s+1})$. Clearly, \mbox{(\ref{eq:ns_HGfixedpointa}--d)} is a linear
system for the unknowns $(\vec U^{m+1,s+1}, P^{m+1,s+1},
\vec X^{m+1,s+1}, \kappa^{m+1,s+1})$ that is of the form
(\ref{eq:ns_algebraic}).

\begin{rem} \label{rem:nondivfree}
At times we will be interested in numerical simulations for
non divergence-free solutions, in particular when considering our convergence
experiments in \S\ref{sec:ns_convergence_results}. Here the condition
{\rm (\ref{eq:ns_mass_bis})} is replaced by
$\nabla\,.\,\vec u = f_{\rm{div}}$ in $\Omega_\pm(t)$.
Clearly, this leads to an additional term on the right-hand side of the
discrete continuity equation {\rm (\ref{eq:ns_HGb})}. In particular,
we have to replace {\rm (\ref{eq:LAb})} with
\begin{equation}\label{eq:LAb_nondivfree}
\left(\nabla\,.\,\vec U^{m+1}, \varphi\right) =
\left(\,f_{\rm{div}},
\varphi-\frac{\left(\varphi,1\right)}{\left(1,1\right)}\right)+
\frac{\left(\varphi, 1\right)}{\vol(\Omega)}\, \int_{\partial_1\Omega}
(\vec I^m_2\,\vec g) \,.\, \unitn \dH{d-1} \quad \forall\ \varphi \in
\pspace^m\,.
\end{equation}
The same adjustment is made on the right-hand side of
{\rm (\ref{eq:ns_HGfixedpointb})}.
In addition, we observe that the weak formulation
{\rm (\ref{eq:ns_weaka_antisym})} needs to be suitably adjusted,
with the additional term $\tfrac{1}{2}\,(\rho\,f_{\rm{div}}\,\vec u,\vec\xi)$
on the right-hand side of {\rm (\ref{eq:ns_weaka_antisym})}, see
\cite{Agnese} for details.
As a consequence, the term
$\frac12\,(\rho^m\,f_{\rm{div}}\,\vec U^m,\vec\xi)$ needs to be added
to the right-hand side of {\rm (\ref{eq:ns_HGa_antisym})} for the scheme
$(\schemeB)$.
\end{rem}

\subsubsection{Solvers for the linear (sub)problems}
Following the authors' approach in \cite{stokesfitted}, see also
\cite{fluidfbp,spurious},
for the solution of (\ref{eq:ns_algebraic}) we use a Schur complement
approach that eliminates $(\kappa^{m+1}, \delta \vec X^{m+1})$ from
(\ref{eq:ns_algebraic}), and then use an iterative solver for the remaining
system in $(\vec U^{m+1}, P^{m+1})$. This approach has the advantage that for
the reduced system well-known solution methods for finite element
discretizations for standard Navier--Stokes discretizations may be employed.
The desired Schur complement can be obtained as follows. Let
\begin{equation} \label{eq:Xi}
\Xi_\Gamma:= \begin{pmatrix}
 0 & - \frac1{\tau}\,\vec N_\Gamma^T \\
\vec N_\Gamma & \vec A_\Gamma
\end{pmatrix} \,,
\end{equation}
be the interface condensed operator. Then (\ref{eq:ns_algebraic}) can be
reduced to
\begin{subequations}
\begin{align}\label{eq:SchurkX}
&
\begin{pmatrix}
\vec B_\Omega + \gamma\,(\Nbulk \ 0)\,\Xi_\Gamma^{-1}\,
\binom{\NbulkT}{0} & \vec C_\Omega \\
\vec C_\Omega^T & 0
\end{pmatrix}
\begin{pmatrix}
\vec U^{m+1} \\ P^{m+1}
\end{pmatrix}
= 
\begin{pmatrix}
\vec c
+\gamma\,(\Nbulk \ 0)\, \Xi_\Gamma^{-1}\,
\binom{0}{-\vec A_\Gamma\,\vec X^m} \\
\vec \beta
\end{pmatrix}
\end{align}
and
\begin{equation}
\binom{\kappa^{m+1}}{\delta\vec X^{m+1}} = \Xi_\Gamma^{-1}\,
\binom{-\NbulkT\,\vec U^{m+1}}{-\vec A_\Gamma\,\vec X^m}\,.
\label{eq:SchurkXb}
\end{equation}
\end{subequations}

We notice that when $\gamma = 0$ the linear system (\ref{eq:SchurkX}) becomes
\begin{equation} \label{eq:ns_system}
\begin{pmatrix}
\vec B_\Omega & \vec C_\Omega \\
\vec C_\Omega^T & 0
\end{pmatrix}
\begin{pmatrix}
\vec U^{m+1} \\ P^{m+1}
\end{pmatrix}
= \begin{pmatrix}
\vec c \\
\vec \beta
\end{pmatrix}\,,
\end{equation}
which for $\mu_+=\mu_-$ and $\rho_+ = \rho_-$
corresponds to a discretization of the one-phase
Navier--Stokes problem. This system is well known and a classical way to solve
it is to use a preconditioned GMRES iteration. Therefore, we solve
(\ref{eq:SchurkX}), which is a slight modification of (\ref{eq:ns_system}),
employing a preconditioned GMRES iterative solver. For the inverse
$\Xi_\Gamma^{-1}$ we employ a sparse $LU$ decomposition, which we obtain with
the help of the sparse factorization package UMFPACK, see \cite{Davis04}, in
order to reuse the factorization. Having obtained $(\vec U^{m+1}, P^{m+1})$ from
(\ref{eq:SchurkX}), we solve (\ref{eq:SchurkXb}) for $(\kappa^{m+1}, \delta\vec
X^{m+1})$.

As preconditioner for (\ref{eq:SchurkX})
we employ the standard Navier--Stokes matrix
\begin{equation}\label{eq:ns_direct_precond}
\begin{pmatrix} \vec B_\Omega & \vec C_\Omega \\
\vec C_\Omega^T & 0 \end{pmatrix},
\end{equation}
recall (\ref{eq:ns_system}).
Of course, as the pressure in the Navier--Stokes problem
(\ref{eq:ns_system}) is only unique up to a constant, the block matrix
(\ref{eq:ns_direct_precond}) is singular. Hence, in order to efficiently apply
the preconditioner (\ref{eq:ns_direct_precond}) as part of an iterative
solution method for (\ref{eq:SchurkX}), one can either use a sparse $QR$
factorization method as provided by SPQR, see \cite{Davis11}. Alternatively,
the matrix (\ref{eq:ns_direct_precond}) can be suitably doctored to make it
invertible, and then a sparse $LU$ factorization can be obtained with the help
of UMFPACK, see \cite{Davis04}, as usual. In practice we pursued the latter
approach, as it proved to be slightly more efficient. We refer to
\cite[\S4.6]{Agnese} for more details.

\subsubsection{Mesh generation and smoothing}\label{sec:smoothing}
Given the initial polyhedral surface $\Gamma^0$, we create a triangulation
$\mathcal{T}^0$ of $\Omega$ that is fitted to $\Gamma^0$ with the help of the
package \verb|Gmsh|, see \cite{GeuzaineR09}.
The mesh generation process is controlled by characteristic length
parameter, which describes the fineness of the resulting mesh.
Then, for $m \geq 0$, having computed the new interface $\Gamma^{m+1}$, we
would like to obtain a bulk triangulation $\mathcal{T}^{m+1}$ that is fitted to
$\Gamma^{m+1}$, and ideally is close to $\mathcal{T}^m$. This is to avoid
unnecessary overhead from remeshing the domain $\Omega$ completely.
To this end, we perform the following smoothing step on $\mathcal{T}^m$, which
is inspired by the method proposed in \cite{Ganesan06}, see also
\cite{GanesanT08}. Having obtained $\delta \vec X^{m+1}$ from the solution of
(\ref{eq:ns_algebraic}), we solve the linear elasticity problem: Find a
displacement $\vec\psi \in [H^1(\Omega)]^d$ such that
\begin{equation}
\nabla\,.\,(2\,\mat D(\vec\psi) + (\nabla\,.\vec\psi)\,\mat\id)
= \vec 0 \quad \mbox{in } \Omega_\pm^m\,,\qquad
\vec\psi = \delta \vec X \quad \mbox{on } \Gamma^m\,, \qquad
\vec\psi\,.\,\unitn = 0 \quad \mbox{on } \partial\Omega\,.
\label{eq:elast}
\end{equation}
In practice we approximate (\ref{eq:elast})
with piecewise linear elements and solve the resulting system of linear
equations with the UMPFACK package, see \cite{Davis04}. The obtained discrete
variant of $\vec\psi$, at every vertex of the current bulk grid
$\mathcal{T}^m$, then represents the variation in their position that we
compute in order to obtain $\mathcal{T}^{m+1}$.

Occasionally the deformation of the mesh becomes too large, and so a
complete remeshing of $\Omega$ becomes necessary. In order to detect the need
for a complete remeshing, we define the angle criterion
\begin{equation}\label{eq:angle_criterion}
\min_{\sigmaO\in {\cal T}^{m+1}}\min_{\alpha\in \measuredangle(\sigmaO)}\alpha
\leq C_a\,,
\end{equation}
where $C_a$ is a fixed constant. Moreover, $\measuredangle(\sigmaO)$ is the set
of all the angles of the simplex $\sigmaO$, which can be computed as
$\alpha_{ij}=\cos^{-1}(-\vec n_i\,.\,\vec n_j)$, $\forall\ i,j\in\{0,\dots,d\}$,
with $\vec n_i$ the unitary normal to the simplex face $i$.
In practice we perform a full remeshing if the
criterion (\ref{eq:angle_criterion}) holds with $C_a=20\degree$.
We stress that in all our numerical simulations we never need to remesh the
interface $\Gamma^m$ itself.

\subsubsection{Velocity interpolation}\label{sec:velocity_interpolation}
We notice that in the three Eulerian schemes $(\schemeAex)$, $(\schemeAim)$
and $(\schemeB)$, all feature at least one term containing the
interpolated velocity $\vec I^m_2\,\vec U^m$. Here the interpolation operator
is needed because the velocity $\vec U^m$, which has been computed at time level
$m-1$, is defined on $\mathcal{T}^{m-1}$. Therefore, at the next time level
$m$, the previous velocity $\vec U^m$ needs to be interpolated on
$\mathcal{T}^m$, which will in general be different from
$\mathcal{T}^{m-1}$, in order to be used.

The interpolation routine simply consists in evaluating the discrete function
$\vec U^m$ in all the degrees of freedom of a piecewise quadratic function
defined on $\mathcal{T}^m$. The efficiency of this computation hinges on
quickly finding the element in the bulk triangulation $\mathcal{T}^{m-1}$
that contains the coordinates of the degree of freedom on $\mathcal{T}^m$
that is currently of interest. Of course, a naive linear looping over all
elements in $\mathcal{T}^{m-1}$ is highly inefficient.
A simple optimization to this naive search consists
in starting a guided search from the last correctly found element, in the hope
that the degrees of freedom on $\mathcal{T}^m$ are traversed in such a way,
the successive points are close to each other. If the element does not
contain the new point, then the guided search moves to the neighbouring
element that is opposite the vertex with the largest distance to the new
point. For a convex domain, this search is always well-defined.
For nonconvex domains this guided search
can be combined with a naive linear search whenever a boundary face is
encountered in the local search for the next neighbour.

\begin{figure}[htbp]
\centering
\includegraphics[width=.25\textwidth]{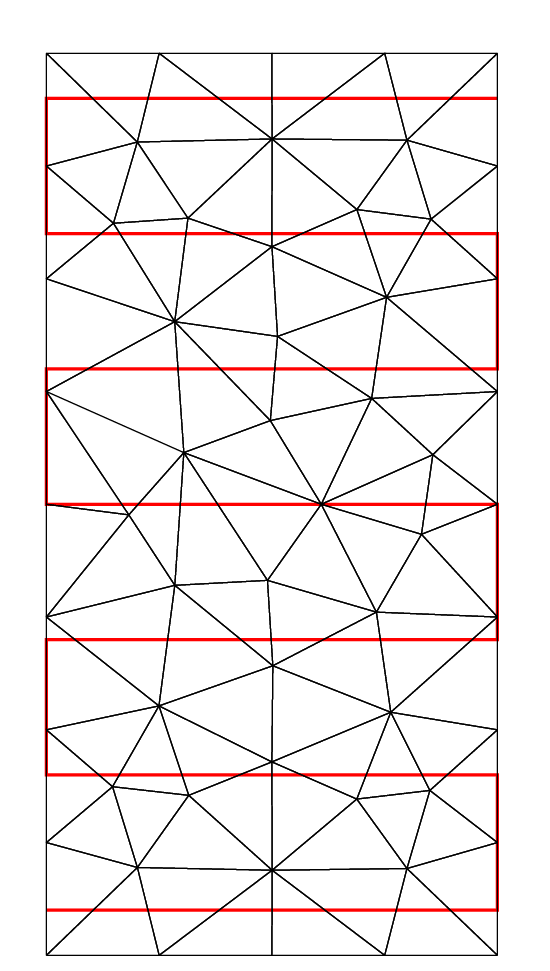}
\caption[Continuous path traversing a fictitious background lattice]
{Continuous path (red) traversing a fictitious background lattice.}
\label{fig:velocity_interpolation_path}
\end{figure}

The last ingredient to our velocity interpolation algorithm is the
traversal of the degrees of freedom on $\mathcal{T}^m$. Ideally we would like
to have a continuous path of elements, which visits, only once, all the
elements of $\mathcal{T}^m$.
This ensures that successive degrees of freedom on the new
mesh are spatially close, and so the local guided search for the next
element of $\mathcal{T}^{m-1}$ to visit is likely to be very short.
Unfortunately, since we use unstructured meshes, such a path, in general,
is difficult to construct. Hence, instead of explicitly constructing such
a path, we create a fictitious background lattice
$\delta\,{\mathbb Z}^d \cap \overline\Omega$ that covers the domain
$\overline\Omega$.
Then each element of $\mathcal{T}^m$ is assigned to the point
on the lattice that is closest to its barycentre.
Depending on the size of the lattice, $\delta$, multiple
entities can be mapped to the same point on the lattice.
Since the lattice is structured, it
is straightforward to build a continuous path which visits all its points. For
example, for a rectangular domain in 2d we can use a snaking path that
traverses the lattice alternating left-to-right and right-to-left, see
Figure~\ref{fig:velocity_interpolation_path}.
Now, in practice, the elements of $\mathcal{T}^m$ are visited following the
given lattice path. This means that, for each lattice point traversed, the
associated entities on $\mathcal{T}^m$ are visited.
Clearly, for the algorithm to be efficient, the lattice size $\delta$ should be
neither too small nor too large.
In practice we use $\delta$ proportional to the characteristic length of the
bulk mesh.

\subsection{ALE scheme} \label{sec:ale_solution_method}
The solution technique used to solve the ALE scheme $(\schemeALE)$ is
very similar to the one used to solve the nonlinear scheme $(\schemeAim)$,
recall (\ref{eq:ns_HGfixedpointa}--d), in that we use a fixed point iteration.
However, the main difference is that at each time step for the ALE scheme we
also find a bulk mesh velocity $\W^{m+1}$, and hence a new bulk mesh
$\Omega^{m+1}$, as part of the solution process.
More precisely, we find a solution for the scheme $(\schemeALE)$ as follows.
Let $\Gamma^m$ and $\vec U^m\in \uspacediscale{g}{m}$ be given.
Set $\vec U^{m+1,0}=\vec U^m$, $\vec \psi^{m+1,0} =
\vec 0$, $\Omega^{m+1,0}=\Omega^m$ and fix $\epsilon_f > 0$.
Then, for $s \geq 0$,
find $(\vec U^{m+1,s+1},P^{m+1,s+1}, \vec X^{m+1,s+1}, \kappa^{m+1,s+1}) \in
\uspacediscale{g}{m}\times \pspaceale^m \times \Vh \times \Wh$ such that
\begin{subequations}
\begin{align}
& \left( \frac{\rho^m}{\tau}\,\vec U^{m+1,s+1}, \vec
\xi\right)_{\Omega^{m+1,s}} + \left(\left(\rho^m\left(\vec
U^{m+1,s}-\frac{\vec\psi^{m+1,s}}{\tau}\right)\,.\,\nabla\right)\,
\vec U^{m+1,s+1}\,, \,\vec \xi \right)_{\Omega^m}\nonumber \\
& \qquad + 2\left(\mu^m\,\mat D(\vec U^{m+1,s+1}), \mat D(\vec \xi)
\right)_{\Omega^m} - \left(P^{m+1,s+1}, \nabla\,.\,\vec \xi\right)_{\Omega^m}
- \gamma\,\left\langle \kappa^{m+1,s+1}\,\vec\nu^m,
\vec\xi\right\rangle_{\Gamma^m} \nonumber \\
& \qquad = \left(\frac{\rho^m}{\tau}\,\vec U^m,\vec \xi
\right)_{\Omega^{m+1,s}}
+ \left(\rho^m\,\vec f_1^{m+1}+\vec f_2^{m+1}, \vec \xi\right)_{\Omega^m}
\quad \forall\ \vec\xi \in \uspacediscale{0}{m}\,, \label{eq:ns_HGalea} \\
& \left(\nabla\,.\,\vec U^{m+1,s+1}, \varphi\right)_{\Omega^m} =
 \frac{\left(\varphi, 1\right)_{\Omega^m}}{\vol(\Omega)}\,
\int_{\partial_1\Omega^m} (\vec I^m_2\,\vec g) \,.\, \unitn \dH{d-1}
\quad \forall\ \varphi \in \pspaceale^m\,,\label{eq:ns_HGaleb} \\
& \left\langle \frac{\vec X^{m+1,s+1} - \vec\id}{\tau} ,\chi\,\vec\nu^m
\right\rangle_{\Gamma^m}^h - \left\langle \vec U^{m+1,s+1}, \chi\,\vec\nu^m
\right\rangle_{\Gamma^m} = 0 \quad \forall\ \chi \in \Wh\,,
\label{eq:ns_HGalec}\\
& \left\langle \kappa^{m+1,s+1}\,\vec\nu^m, \vec\eta \right\rangle_{\Gamma^m}^h
+ \left\langle \nabs\,\vec X^{m+1,s+1}, \nabs\,\vec \eta
\right\rangle_{\Gamma^m} = 0 \quad \forall\ \vec\eta \in \Vh\,.
\label{eq:ns_HGaled}
\end{align}
\end{subequations}
Having obtained $\vec X^{m+1,s+1}$, the new bulk displacement
$\vec \psi^{m+1,s+1}$ is computed via the linear elasticity problem
(\ref{eq:elast}), see \S\ref{sec:smoothing},
and the new bulk mesh $\Omega^{m+1,s+1}$ is updated accordingly.
This iteration is repeated until $\|U^{m+1,s+1}-U^{m+1,s}\|_{L^\infty}
\leq\epsilon_f$ and $\|\vec \psi^{m+1,s+1}-\vec \psi^{m+1,s}\|_{L^\infty}
\leq\epsilon_f$.
Finally, set
$\Gamma^{m+1} = \vec X^{m+1,s+1}(\Gamma^m)$, $\Omega^{m+1,s+1}=\Omega^{m+1}$
and \linebreak $(\vec U^{m+1}, P^{m+1}, \kappa^{m+1},\W^{m+1}) =
(\vec U^{m+1,s+1}, P^{m+1,s+1},\kappa^{m+1,s+1},\tfrac1\tau\,\psi^{m+1,s+1})$.
Clearly, \linebreak
\mbox{(\ref{eq:ns_HGalea}--d)} is once again
a linear system of equations for the unknowns
$(\vec U^{m+1,s+1}, P^{m+1,s+1},$ $ \vec X^{m+1,s+1},
\kappa^{m+1,s+1})$ of the form (\ref{eq:ns_algebraic}).

Similarly to the Eulerian schemes, we monitor the mesh quality of
the bulk meshes $\Omega^{m+1}$ with the help of the criterion
(\ref{eq:angle_criterion}). In practice, if (\ref{eq:angle_criterion})
holds with $C_a=20\degree$, then we perform a full remeshing of the domain
$\Omega$, keeping $\Gamma^{m+1}$ fixed, and obtain the discrete velocity
approximation on the new mesh via the procedure described in
\S\ref{sec:velocity_interpolation}.

\setcounter{equation}{0}
\section{Numerical results} \label{sec:nr}

We implemented our numerical schemes within the \verb|DUNE| framework,
using the FEM toolbox \verb|dune-fem|, see
\cite{dunegridpaperI08,dunegridpaperII08,dunefempaper10}.
As grid manager we employed \verb|ALBERTA|, \cite{Alberta},
and all the meshes were created with the library \verb|Gmsh|,
\cite{GeuzaineR09}. All the schemes described in this paper are implemented as
a \verb|DUNE| module called \verb|dune-navier-stokes-two-phase|, available at
\myurl{https://github.com/magnese/dune-navier-stokes-two-phase}.
The simulations were performed as single thread computations on a cluster of
Xeon CPUs with frequency between 2.40GHz to 3.00GHz, and with at
least 16GB of memory.

\subsection{Convergence experiments}\label{sec:ns_convergence_results}
Here we consider convergence experiments for exact solutions to the
incompressible two-phase Navier--Stokes flow problem. These exact solutions
will be characterised by expanding spheres. In particular, let
$\Gamma(t) = \{ \vec z \in \R^d : |\vec z\,| = r(t)\}$ be a sphere of radius
$r(t)$, with curvature $\varkappa(t) = -\,\frac{d-1}{r(t)}$.
Then, for $\alpha \in \R_{\geq0}$, it can be shown that the expanding sphere
with radius
\begin{subequations}
\begin{equation} \label{eq:ns_sol1a}
r(t) = \mathrm{e}^{\alpha\,t}\,r(0)\,,
\end{equation}
together with
\begin{equation} \label{eq:ns_sol1b}
\vec u(\vec z, t) = \alpha\,\vec z \,, \quad
p(\vec z, t) = -\bigg[\gamma - 2\,\alpha\,\frac{\mu_+ - \mu_-}
{d-1}r(t)\bigg]\,\varkappa(t)\left[ \charfcn{\Omega_-(t)} -
\frac{\vol(\Omega_-(t))}{\vol(\Omega)}\right],
\end{equation}
\end{subequations}
is an exact solution to the problem (\ref{eq:ns_full_momentum}--g), with
(\ref{eq:ns_full_mass}) replaced by
\[
\nabla\,.\,\vec u = \alpha\,d \quad \mbox{in } \Omega_\pm(t)\,,
\]
and with $\vec f(\cdot, t) = \rho\,\alpha^2\,\vec\id$ and
$\vec g = \alpha\,\vec\id$ on e.g.\ $\Omega = (-1,1)^d$ with
$\partial_1\Omega=\partial\Omega$.

A nontrivial divergence free and radially symmetric solution $\vec u$
can be constructed on a domain that does not contain the origin. To this end,
consider e.g.\ $\Omega = (-1,1)^d \setminus [-\frac13, \frac13]^d$.
Then, for $\alpha \in \R_{\geq0}$, the expanding sphere with radius
\begin{subequations}
\begin{equation} \label{eq:ns_sol2a}
r(t) = ([r(0)]^d + \alpha\,t\,d)^\frac1d \,,
\end{equation}
together with
\begin{equation} \label{eq:ns_sol2b}
\vec u(\vec z, t) = \alpha\,\frac{\vec z}{|\vec z\,|^d}\,, \quad
p(\vec z, t) = -\,\bigg(\gamma +2\,\alpha\,\frac{\mu_+ - \mu_-}
{r(t)^{d-1}}\bigg)\,\varkappa(t)\left[ \charfcn{\Omega_-(t)} -
\frac{\vol(\Omega_-(t))}{\vol(\Omega)}\right],
\end{equation}
\end{subequations}
is an exact solution to the problem (\ref{eq:ns_full_momentum}--g) with
$\vec f(\vec z, t) = \rho\,(1-d)\,\alpha^2\,|\vec z|^{-2d}\,\vec z$ and
$\vec g(\vec z) = \alpha\,|\vec z|^{-d}\,\vec z$ on
$\partial_1\Omega=\partial\Omega$.

For the following convergence tests,
we choose the initial surface $\Gamma(0) = \{ \vec z \in \R^d : |\vec z| =
\frac12 \}$, $\partial_1\Omega=\partial\Omega$ and
adaptive bulk meshes that use a finer resolution close to the interface. We
compute the discrete solution over the time interval $[0,1]$.
We begin with the exact solution (\ref{eq:ns_sol1a},b) for $\alpha = 0.15$
on $\Omega = (-1,1)^2$, with the physical parameters
$\rho_+ = \rho_- = \mu_+ = \mu_- = \gamma = 1$.
Details on the discretization parameters are given in
Table~\ref{tab:nsexpandingbubbleelements}, where we also state the final
number of bulk elements, $J_\Omega^M$, for the scheme $(\schemeAim)$.
The values of $J_\Omega^M$ for the remaining three schemes were very similar.
\begin{table*}
\center
\begin{tabular}{rrrrrr}
\hline
& & \multicolumn{2}{c}{(\ref{eq:ns_sol1a},b)} &
    \multicolumn{2}{c}{(\ref{eq:ns_sol2a},b)} \\
$J_\Gamma$ & $\tau$
& $J_\Omega^0$ & $J_\Omega^M$ & $J_\Omega^0$ & $J_\Omega^M$ \\
\hline
 32 & $6.4\cdot10^{-2}$ &     296 &     310 &  460   &  216 \\
 64 & $1.6\cdot10^{-2}$ &  1\,240 &  1\,240 & 1\,040 &  444 \\
128 &   $4\cdot10^{-3}$ &  4\,836 &  4\,836 & 2\,628 & 1\,378 \\
256 &         $10^{-3}$ & 18\,476 & 18\,476 & 7\,460 & 4\,476 \\
\hline
\end{tabular}
\caption[Navier--Stokes expanding bubble meshes parameters]
{Discretization parameters for the convergence experiments. The values for
$J_\Omega^M$ correspond to the scheme $(\schemeAim)$.}
\label{tab:nsexpandingbubbleelements}
\end{table*}%
We report on the computed errors in Table~\ref{tab:nsexpandingbubbleIp2p1p0},
where we have defined the error quantities
$\errorXx := \max_{m=1,\ldots, M} \|\Gamma^m - \Gamma(t_m)\|_{L^\infty}$,
where $\|\Gamma^m - \Gamma(t_m)\|_{L^\infty} :=
\max_{k=1,\ldots, K_\Gamma} {\rm dist}( \vec q^m_k, \Gamma(t_m))$,
as well as
$\LerrorUu2 := \left[\tau\,\sum_{m=1}^M \|\vec U^m - I^m_2\,\vec u(\cdot,
t_m)\|_{L^2(\Omega)}^2 \right]^\frac12$,
$\HerrorUu2 := $ \linebreak
$\left[\tau\,\sum_{m=1}^M \|\vec U^m - I^m_2\,\vec u(\cdot,
t_m)\|_{H^1(\Omega)}^2 \right]^\frac12$
and
$\LerrorPp := \left[\tau\,\sum_{m=1}^M \|P^m - p(\cdot,t_m)\|_{L^2(\Omega)}^2
\right]^\frac12$.
Moreover, we also define the estimated order of convergence EOC as
$\mbox{EOC}={\ln{\frac{\mbox{error}_1}{\mbox{error}_0}}}/
{\ln{\frac{\mbox{h}_1}{\mbox{h}_0}}}$,
where $\mbox{error}_1$ and $\mbox{error}_0$ are the errors for the computations
with characteristic lengths $\mbox{h}_1 < \mbox{h}_0$.
\begin{table*}
\center
\hspace*{-3.25cm}
\begin{tabular}{rllllllr}
\hline
$J_\Gamma$ & $\errorXx$ & $\LerrorUu2$ & EOC & $\HerrorUu2$ & $\LerrorPp$ & EOC
& CPU[s] \\
\hline
& \multicolumn{7}{c}{$(\schemeAex)$} \\
\hline
 32 & 3.96456e-04 & 0 & - & 0 & 3.05157e-01 &    - &     5 \\
 64 & 1.03429e-04 & 0 & - & 0 & 1.57053e-01 & 0.96 &    61 \\
128 & 2.61302e-05 & 0 & - & 0 & 7.09596e-02 & 1.15 &  1\,317 \\
256 & 6.53463e-06 & 0 & - & 0 & 1.99794e-02 & 1.79 & 29\,194 \\
\hline
& \multicolumn{7}{c}{$(\schemeAim)$} \\
\hline
 32 & 3.96456e-04 & 0 & - & 0 & 3.05157e-01 &    - &     6 \\
 64 & 1.03429e-04 & 0 & - & 0 & 1.57053e-01 & 0.96 &   134 \\
128 & 2.61302e-05 & 0 & - & 0 & 7.09596e-02 & 1.15 &  2\,649 \\
256 & 6.53463e-06 & 0 & - & 0 & 1.99794e-02 & 1.79 & 49\,990 \\
\hline
& \multicolumn{7}{c}{$(\schemeB)$} \\
\hline
 32 & 3.96456e-04 & 0 & - & 0 & 3.05157e-01 &    - &       5 \\
 64 & 1.03429e-04 & 0 & - & 0 & 1.57053e-01 & 0.96 &      45 \\
128 & 2.61302e-05 & 0 & - & 0 & 7.09596e-02 & 1.15 &  1\,567 \\
256 & 6.53463e-06 & 0 & - & 0 & 1.99794e-02 & 1.79 & 30\,997 \\
\hline
& \multicolumn{7}{c}{$(\schemeALE)$} \\
\hline
 32 & 3.96823e-04 & 1.94467e-06 &    - & 1.74315e-05 & 3.01900e-01 &    - &
10 \\
 64 & 1.03462e-04 & 1.27547e-07 & 3.93 & 1.33117e-06 & 1.57054e-01 & 0.94 &
160 \\
128 & 2.61390e-05 & 4.26296e-08 & 1.58 & 3.78211e-07 & 7.09597e-02 & 1.15 &
2\,420 \\
256 & 6.53680e-06 & 1.09827e-08 & 1.91 & 9.46271e-08 & 1.99793e-02 & 1.79 &
47\,181 \\
\hline
\end{tabular}
\hspace*{-3.25cm}
\caption[Navier--Stokes expanding bubble I errors]
{($\rho_+ = \rho_- = \mu_+ = \mu_- = \gamma = 1,\alpha=0.15$)
Convergence experiments for (\ref{eq:ns_sol1a},b)
on $(-1,1)^2$ over the time interval $[0,1]$, with the
discretization parameters from Table~\ref{tab:nsexpandingbubbleelements}.}
\label{tab:nsexpandingbubbleIp2p1p0}
\end{table*}%
We observe from Table~\ref{tab:nsexpandingbubbleIp2p1p0}
that the schemes $(\schemeAex)$, $(\schemeAim)$ and $(\schemeB)$ capture the
velocity solution exactly. This is possible because the true velocity is linear,
see (\ref{eq:ns_sol1b}). We also notice that all the schemes capture the
interface position with the same accuracy and exhibit similar errors for the
pressure. In terms of the overall CPU times, we notice that the two linear
schemes $(\schemeAex)$ and $(\schemeB)$, as well as the two nonlinear schemes
$(\schemeAim)$ and $(\schemeALE)$, have similar execution times, with the CPU
times for the linear schemes lower than for the nonlinear schemes. Finally, we
note that the number of full remeshings performed when $J_\Gamma = 32$ was,
depending on the scheme, between 1 and 4  while it was always 0 in all the
other simulations for this test case.

In a second set of convergence experiments, we fix
$\Omega = (-1,1)^2 \setminus[-\frac13,\frac13]^2$
for the exact solution (\ref{eq:ns_sol2a},b) with $\alpha=0.15$ and
the physical parameters
\[
\rho_+ = 10^3\,,\quad \rho_- = 10^2\,,\quad \mu_+ = 10\,,\quad \mu_- = 1\,,\quad
\gamma = 1\,.
\]
We report on the errors of the exact solution (\ref{eq:ns_sol2a},b) in
Table~\ref{tab:nsexpandingbubbleIIp2p1p0}.
\begin{table*}
\center
\hspace*{-3.25cm}
\begin{tabular}{rllllllr}
\hline
$J_\Gamma$ & $\errorXx$ & $\LerrorUu2$ & EOC & $\HerrorUu2$ & $\LerrorPp$ & EOC
& CPU[s] \\
\hline
& \multicolumn{7}{c}{$(\schemeAex)$} \\
\hline
 32 & 4.13976e-03 & 1.24661e-03 &    - & 2.59441e-02 & 2.28403e-00 &    - &
8 \\
 64 & 1.07627e-03 & 4.80240e-04 & 1.38 & 1.35253e-02 & 1.20439e-00 & 0.92 &
102 \\
128 & 2.55529e-04 & 3.70025e-04 & 0.38 & 1.20309e-02 & 5.89258e-01 & 1.03 &
2\,810 \\
256 & 6.66480e-05 & 1.42910e-04 & 1.34 & 6.48222e-03 & 2.69953e-01 & 1.10 &
88\,056 \\
\hline
& \multicolumn{7}{c}{$(\schemeAim)$} \\
\hline
 32 & 3.99899e-03 & 1.23928e-03 &    - & 2.61424e-02 & 2.30983e-00 &    - &
11 \\
 64 & 1.07657e-03 & 4.80445e-04 & 1.37 & 1.35256e-02 & 1.20449e-00 & 0.94 &
126 \\
128 & 2.55562e-04 & 3.70165e-04 & 0.38 & 1.20325e-02 & 5.89236e-01 & 1.03 &
3\,223 \\
256 & 6.66478e-05 & 1.42914e-04 & 1.34 & 6.48232e-03 & 2.69982e-01 & 1.10 &
95\,315 \\
\hline
& \multicolumn{7}{c}{$(\schemeB)$} \\
\hline
 32 & 5.56564e-02 & 3.33673e-01 &    - & 8.38557e-00 & 9.98404e+01 &    - &
8 \\
 64 & 1.05282e-03 & 4.75502e-04 & 9.45 & 1.37127e-02 & 2.01792e+01 & 2.31 &
112 \\
128 & 2.52742e-04 & 3.82170e-04 & 0.32 & 1.24067e-02 & 2.00690e+01 & 0.01 &
3\,138 \\
256 & 6.55951e-05 & 1.45885e-04 & 1.36 & 6.67192e-03 & 2.00294e+01 &    0 &
98\,893 \\
\hline
& \multicolumn{7}{c}{$(\schemeALE)$} \\
\hline
 32 & 4.20717e-03 & 1.33379e-03 &    - & 3.03240e-02 & 4.92318e-00 &    - &
15 \\
 64 & 1.11280e-03 & 4.91358e-04 & 1.44 & 1.42754e-02 & 2.47229e-00 & 0.99 &
90 \\
128 & 2.54196e-04 & 3.30194e-04 & 0.57 & 1.31294e-02 & 1.24058e-00 & 0.99 &
991 \\
256 & 6.37496e-05 & 1.35898e-04 & 1.25 & 6.41819e-03 & 6.08051e-01 & 1.01 &
11\,970 \\
\hline
\end{tabular}
\hspace*{-3.25cm}
\caption[Navier--Stokes expanding bubble II errors]
{($\rho_+ = 10^3,\rho_- = 10^2,\mu_+ = 10,\mu_- =1,\gamma = 1,\alpha=0.15$)
Convergence experiments for (\ref{eq:ns_sol2a},b)
on $(-1,1)^2\setminus[-\frac{1}{3},\frac{1}{3}]^2$
over the time interval $[0,1]$, with the discretization parameters from
Table~\ref{tab:nsexpandingbubbleelements}.}
\label{tab:nsexpandingbubbleIIp2p1p0}
\end{table*}%
Firstly, we observe that the velocity and the interface position is captured
with similar accuracy by all four schemes. However, we notice that the
pressure error does not appear to converge for the scheme $(\schemeB)$,
which may be caused by the discontinuous jumps in density and viscosity.
In addition, the scheme $(\schemeALE)$ exhibits slightly higher pressure
errors compared to $(\schemeAex)$ and $(\schemeAim)$.
We also note the vastly superior CPU times for the scheme $(\schemeALE)$
compared to the other schemes. This is caused by the fact that
the domain is nonconvex, and so we employ a simple linear search in
our velocity interpolation routine, recall \S\ref{sec:velocity_interpolation}.
Of course, the ALE scheme only needs to interpolate the velocity when a
complete bulk remeshing is performed, and this leads to a significant reduction
in the total simulation time. Finally, we note that the number of full
remeshings performed during any simulation for this test case was always between
3 and 9.

\subsection{Benchmark computations} \label{sec:benchmark_computations}

Here we will present benchmark computations for the two rising bubble
test cases proposed in \cite[Table~I]{HysingTKPBGT09}.
To this end, we use the setup described in
\cite[Figure~2]{HysingTKPBGT09}, which is
$\Omega = (0,1) \times (0,2)$ with $\partial_1\Omega = [0,1] \times \{0,2\}$
and $\partial_2\Omega = \{0,1\} \times (0,2)$. Moreover, the initial
interface is $\Gamma_0 = \{\vec z \in \R^2 : |\vec z - (\tfrac{1}{2},
\tfrac{1}{2})^T| = \frac{1}{4}\}$. We also let
$\vec f = -0.98\,\rho\,\vec\ek_2$ and $\vec g=\vec 0$.
For the discretization parameters, we choose
$\tau=10^{-3}$, $T=3$ and nearly uniform bulk meshes.

We recall that the paper \cite{HysingTKPBGT09} used the proposed benchmarks to
compare two different codes based on the level-set method with a direct front
tracking ALE method, while the same benchmarks were used in
\cite{AlandV12,GarckeHK16} to
investigate different phase-field methods. Finally,
in \cite{fluidfbp} an unfitted
front tracking method was used for these benchmark problems.

\subsubsection{Benchmark problem I}
The physical parameters for the test case 1 in \cite[Table~I]{HysingTKPBGT09}
are given by
\begin{equation} \label{eq:Hysing1}
\rho_+ = 10^3\,,\quad \rho_- = 10^2\,,\quad \mu_+ = 10\,,\quad \mu_- = 1\,,\quad
\gamma = 24.5\,.
\end{equation}
We list our discretization parameters for this benchmark computation
in Table~\ref{tab:risingbubble2Delements}, and report on the quantitative
results for our four schemes in Table~\ref{tab:risingbubbleIp2p1p0}.
Here we have defined the $\vec\ek_d$-component of the bubble's centre of mass
as $z_c^m = \int_{\Omega_-^m} \vec\id\,.\,\vec\ek_d \dL{d} / \vol(\Omega_-^m)$,
the degree of sphericity as
$\strikes^m ={\pi^{\frac{1}{d}}[2\,d\,\vol(\Omega_-^m)]^\frac{d-1}{d}}/
{\surfvol(\Gamma^m)}$, and the bubble's rise velocity as
$V^m_c = \int_{\Omega_-^m}\vec U^m\,.\,\vec\ek_d \dL{d} / \vol(\Omega_-^m)$.
\begin{table*}
\center
\begin{tabular}{rrrrrr}
\hline
$J_\Gamma$ & $J_\Omega^0$ & \multicolumn{4}{c}{$J_\Omega^M$} \\ \hline
& & $(\schemeAex)$ & $(\schemeAim)$ & $(\schemeB)$ & $(\schemeALE)$ \\
\hline
 32 &  2\,210 &  1\,816 &  1\,816 &  1\,794 &  1\,820 \\
 64 &  8\,822 &  7\,544 &  7\,156 &  7\,490 &  7\,566 \\
128 & 35\,092 & 29\,276 & 29\,014 & 28\,968 & 28\,296 \\
\hline
\end{tabular}
\caption{Discretization parameters for the benchmark problem I.}
\label{tab:risingbubble2Delements}
\end{table*}%
\begin{table*}
\center
\hspace*{-3.25cm}
\begin{tabular}{r|rrr|rrr|}
\hline
 & $J_\Gamma=32$ & $J_\Gamma=64$ & $J_\Gamma=128$
 & $J_\Gamma=32$ & $J_\Gamma=64$ & $J_\Gamma=128$ \\ \hline
& \multicolumn{3}{c|}{$(\schemeAex)$} & \multicolumn{3}{c|}{$(\schemeAim)$} \\
\cmidrule{2-7}
$\strikes_{\min}$                & 0.8929 & 0.8975 & 0.9001  & 0.8925 & 0.8973 & 0.9004 \\
$t_{\strikes = \strikes_{\min}}$ & 1.9040 & 1.9040 & 1.9170  & 2.0410 & 1.9120 & 1.9160 \\
$V_{c,\max}$                     & 0.2439 & 0.2424 & 0.2411  & 0.2439 & 0.2423 & 0.2410 \\
$t_{V_c = V_{c,\max}}$           & 0.9350 & 0.9300 & 0.9180  & 0.9340 & 0.9300 & 0.9190 \\
$z_c(t=3)$                       & 1.0829 & 1.0852 & 1.0822  & 1.0820 & 1.0841 & 1.0819 \\
CPU[s]                           &  10\,236 &  29\,610 & 234\,396  &  17\,771 &  72\,837 & 403\,904 \\
\cmidrule{2-7}
& \multicolumn{3}{c|}{$(\schemeB)$} & \multicolumn{3}{c|}{$(\schemeALE)$} \\
\cmidrule{2-7}
$\strikes_{\min}$                & 0.8788 & 0.8851 & 0.8883  & 0.9013 & 0.9010 & 0.9030 \\
$t_{\strikes = \strikes_{\min}}$ & 1.7920 & 1.7190 & 1.6790  & 1.9460 & 1.9050 & 1.8460 \\
$V_{c,\max}$                     & 0.2359 & 0.2344 & 0.2334  & 0.2435 & 0.2422 & 0.2410 \\
$t_{V_c = V_{c,\max}}$           & 0.8860 & 0.8820 & 0.8750  & 0.9370 & 0.9330 & 0.9200 \\
$z_c(t=3)$                       & 1.0648 & 1.0622 & 1.0607  & 1.0877 & 1.0865 & 1.0824 \\
CPU[s]                           &   5\,887 &  19\,703 & 146\,538  &  13\,628 &  63\,136 & 351\,867 \\
\hline
\end{tabular}
\hspace*{-3.25cm}
\caption{($\rho_+ = 10^3,\rho_- = 10^2,\mu_+ = 10,\mu_- =1,\gamma = 24.5$)
Benchmark problem I, with the discretization parameters from
Table~\ref{tab:risingbubble2Delements}.}
\label{tab:risingbubbleIp2p1p0}
\end{table*}%
We observe that the results in Table~\ref{tab:risingbubbleIp2p1p0} are in very
good agreement, except for the scheme $(\schemeB)$. For this reason, we do
not consider the scheme $(\schemeB)$ further in this section.
The results for the remaining three schemes agree very well with
the corresponding
numbers from the finest discretization run of group 3 in \cite{HysingTKPBGT09},
which are given by 0.9013, 1.9000, 0.2417, 0.9239 and 1.0817,
and which are generally accepted as the most accurate in that paper.
In fact, very similar results are also obtained in
\cite[Tables ~2 and 3]{fluidfbp}.
Moreover, we note that the number of full remeshings performed during
any simulation for this test case was always between 2 and 4.

\begin{figure}[htbp]
\centering
\subfloat[$t=10^{-1}$]{\includegraphics[width=.45\textwidth]
{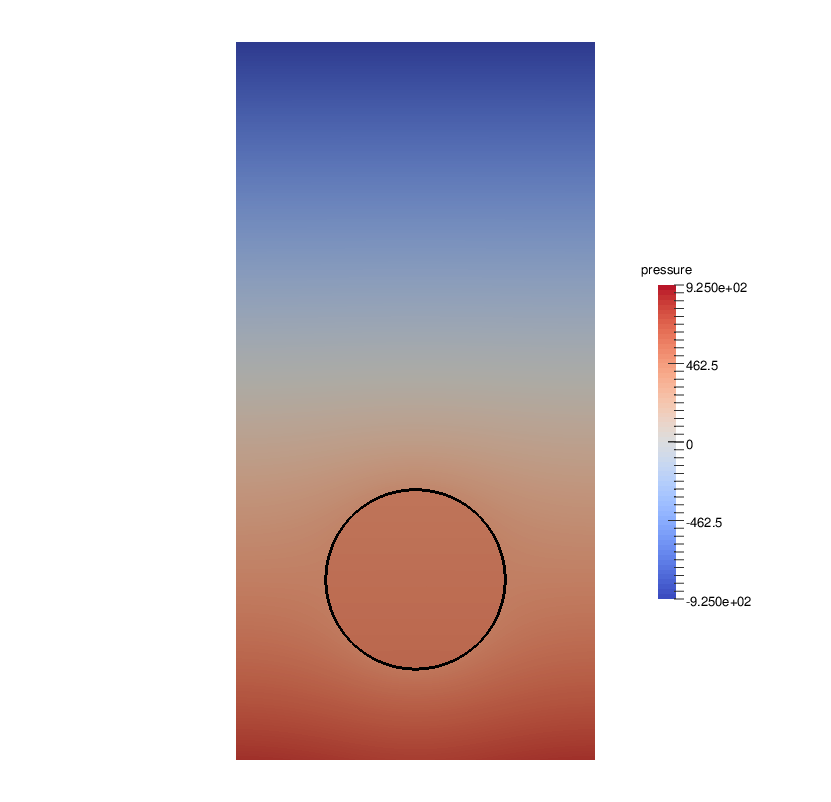}}
\subfloat[$t=1$]{\includegraphics[width=.45\textwidth]
{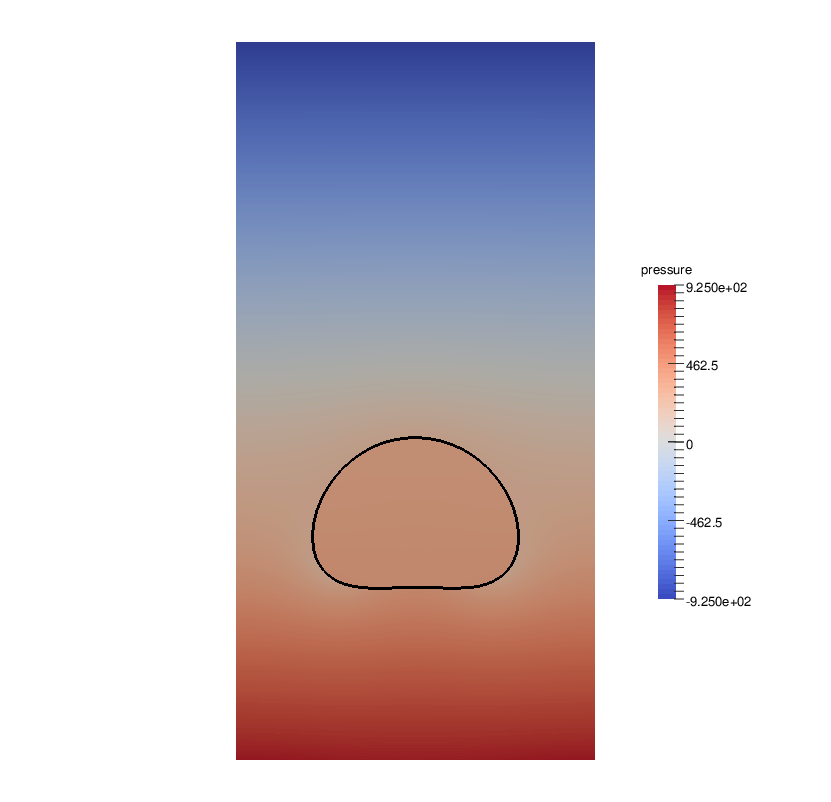}}\\
\subfloat[$t=2$]{\includegraphics[width=.45\textwidth]
{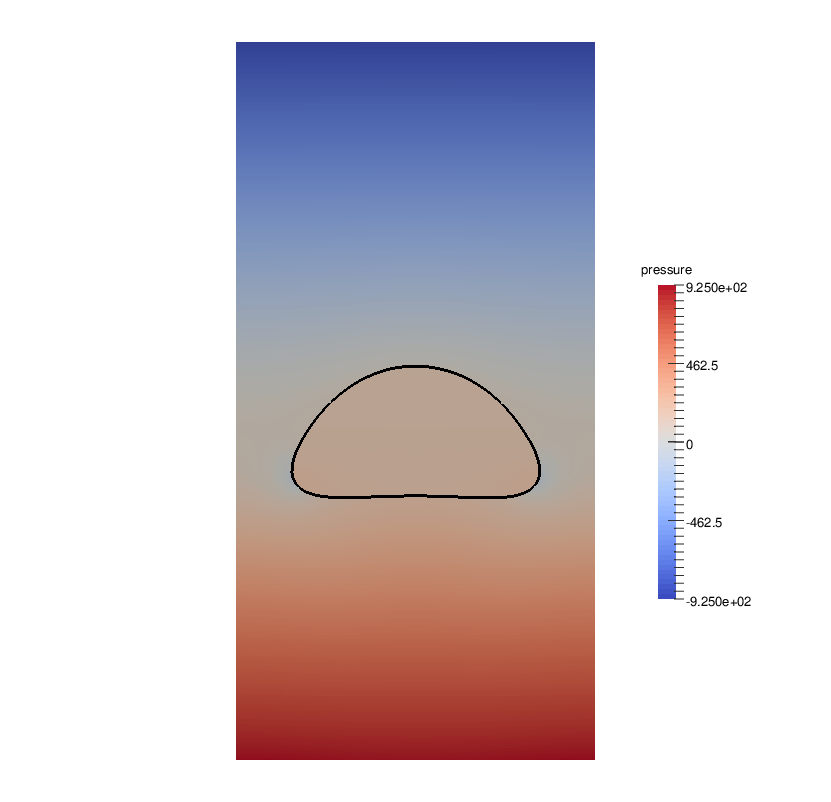}}
\subfloat[$t=3$]{\includegraphics[width=.45\textwidth]
{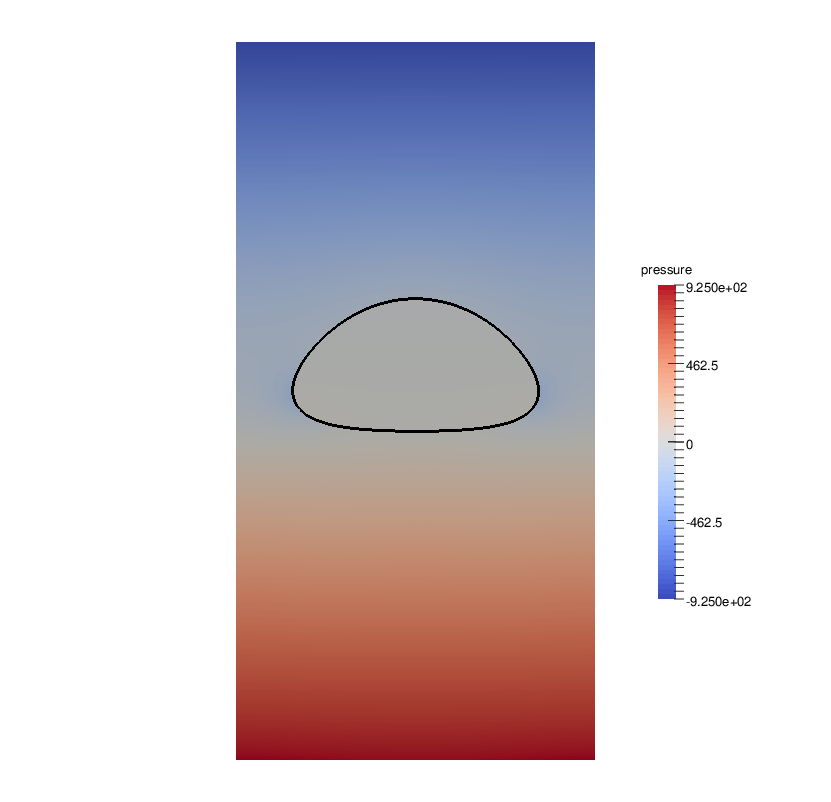}}
\caption{($\rho_+ = 10^3,\rho_- = 10^2,\mu_+ = 10,\mu_- =1,\gamma = 24.5$)
Pressure evolution for the benchmark problem I for the scheme $(\schemeAim)$
with $J_\Gamma=128$ interface elements.}
\label{fig:risingbubbleIpressure}
\end{figure}

\begin{figure}[htbp]
\centering
\subfloat[$t=10^{-1}$]{\includegraphics[width=.45\textwidth]
{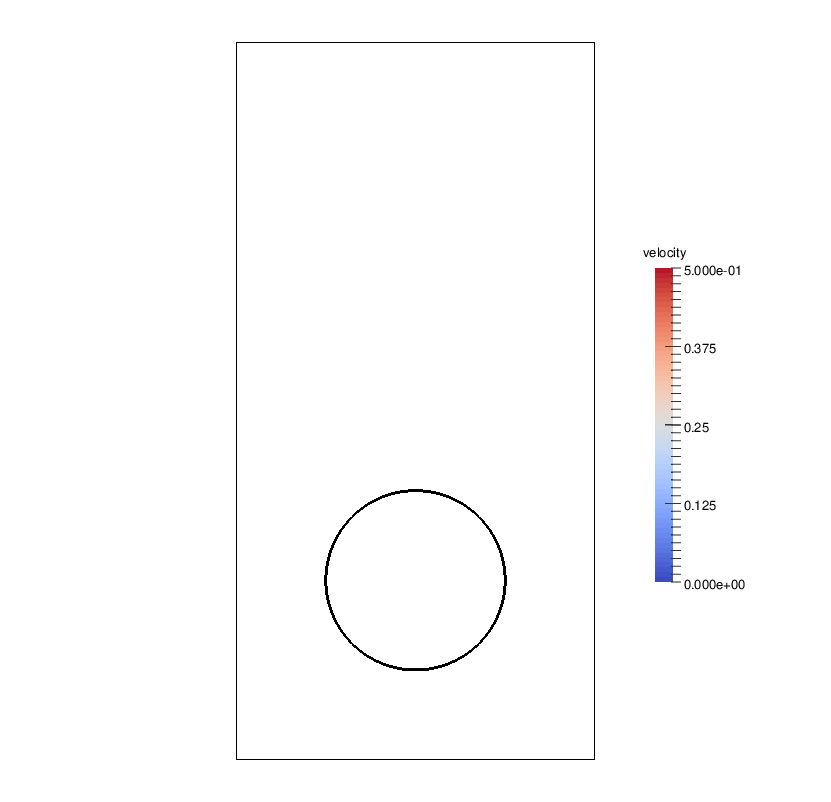}}
\subfloat[$t=1$]{\includegraphics[width=.45\textwidth]
{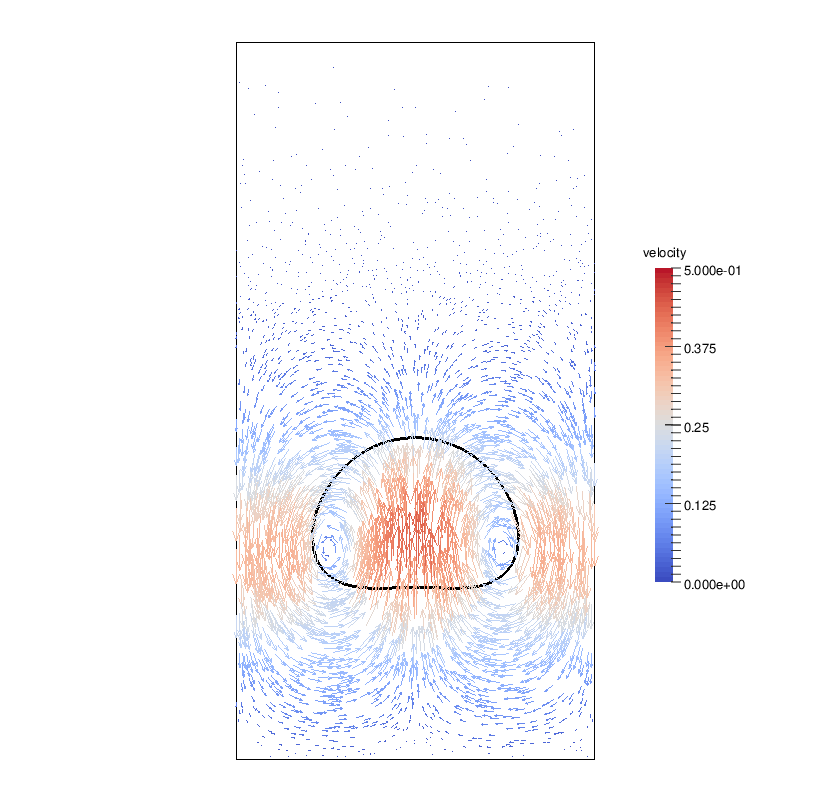}}\\
\subfloat[$t=2$]{\includegraphics[width=.45\textwidth]
{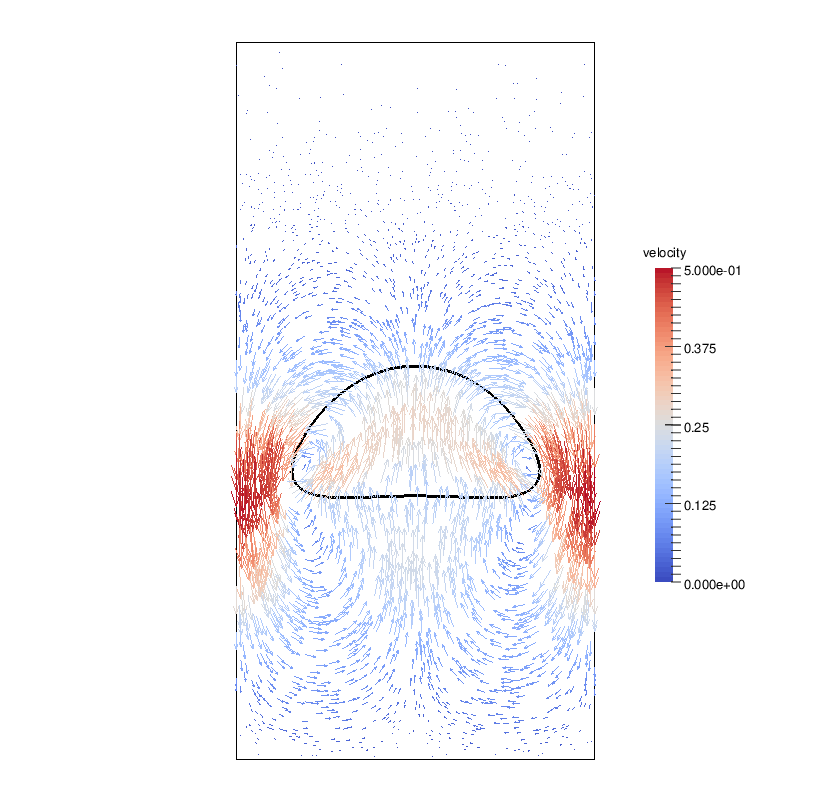}}
\subfloat[$t=3$]{\includegraphics[width=.45\textwidth]
{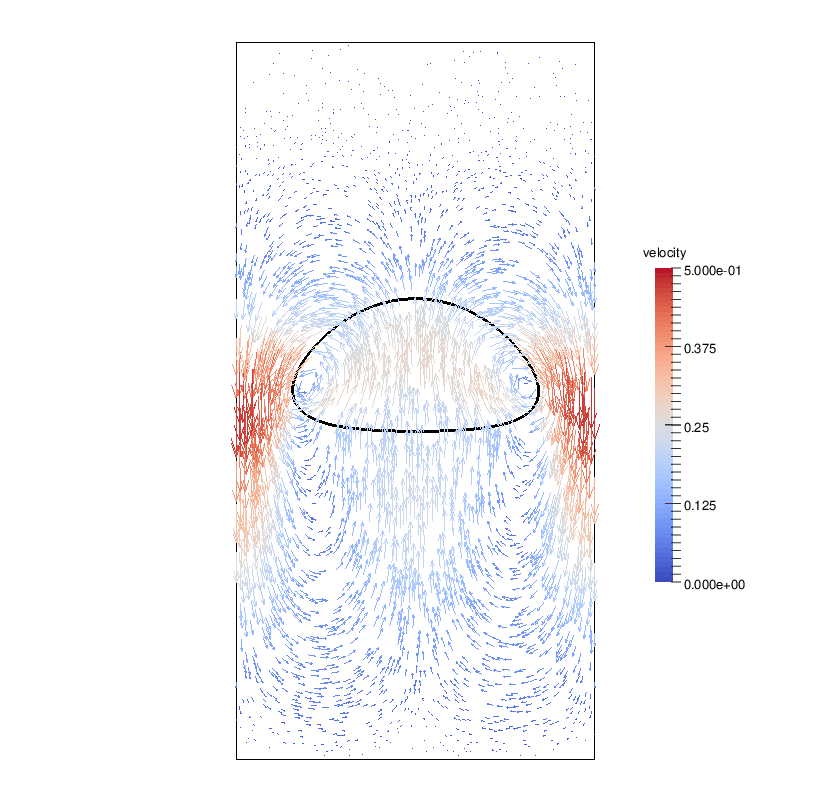}}
\caption{Velocity vector field for the simulation in
Figure~\ref{fig:risingbubbleIpressure}.}
\label{fig:risingbubbleIvelocity}
\end{figure}

\begin{figure}[htbp]
\centering
\includegraphics[width=.45\textwidth]
{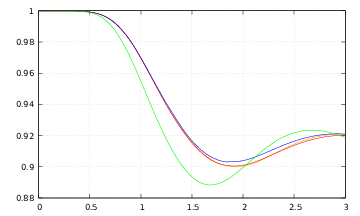}
\caption{A plot of the sphericity $\strikes$ over time for the simulation in
Figure~\ref{fig:risingbubbleIpressure} for the schemes
$(\schemeAex)$ (orange), $(\schemeAim)$ (red), $(\schemeB)$ (green) and
$(\schemeALE)$ (blue).}
\label{fig:risingbubbleIsphericity}
\end{figure}

\begin{figure}[htbp]
\centering
\includegraphics[width=.45\textwidth]
{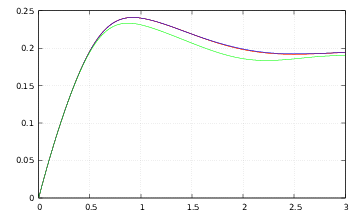}
\caption{A plot of the rising velocity $V_c$ over time for the simulation in
Figure~\ref{fig:risingbubbleIpressure} for the schemes
$(\schemeAex)$ (orange), $(\schemeAim)$ (red), $(\schemeB)$ (green) and
$(\schemeALE)$ (blue).}
\label{fig:risingbubbleIrisingvelocity}
\end{figure}

\begin{figure}[htbp]
\centering
\includegraphics[width=.45\textwidth]
{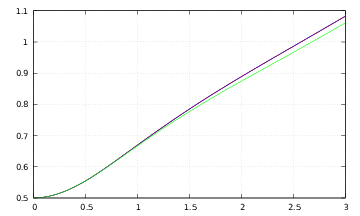}
\caption{A plot of the barycentre $z_c$ over time for the simulation in
Figure~\ref{fig:risingbubbleIpressure} for the schemes
$(\schemeAex)$ (orange), $(\schemeAim)$ (red), $(\schemeB)$ (green) and
$(\schemeALE)$ (blue).}
\label{fig:risingbubbleIbarycenter}
\end{figure}

\begin{figure}[htbp]
\centering
\includegraphics[width=.45\textwidth]
{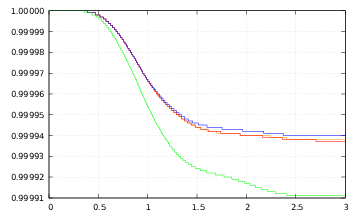}
\caption{A plot of the relative inner area
$\frac{\mathcal{L}^2(\Omega^m_-)}{\mathcal{L}^2(\Omega^0_-)}$ over time for the
simulation in Figure~\ref{fig:risingbubbleIpressure} for the schemes
$(\schemeAex)$ (orange), $(\schemeAim)$ (red), $(\schemeB)$ (green) and
$(\schemeALE)$ (blue).}
\label{fig:risingbubbleIinnervolume}
\end{figure}

In Figure~\ref{fig:risingbubbleIpressure} we show the evolution of the discrete
pressures for a simulation with $J_\Gamma=128$ interface elements using the
scheme $(\schemeAim)$, while the velocities are visualized in
Figure~\ref{fig:risingbubbleIvelocity}.
Finally, in Figures~\ref{fig:risingbubbleIsphericity},
\ref{fig:risingbubbleIrisingvelocity}, and \ref{fig:risingbubbleIbarycenter}
are reported, respectively, the evolution of
the sphericity $\strikes$, rising velocity $V_c$ and barycentre $z_c$ over
time. We notice almost identical results for the three schemes
$(\schemeAex)$, $(\schemeAim)$ and $(\schemeALE)$, which are in very good
agreement with the plots shown in \cite{HysingTKPBGT09}. We also note the
excellent area conservation shown in \ref{fig:risingbubbleIinnervolume},
where we visualize the evolution of the relative
inner area $\frac{\mathcal{L}^2(\Omega^m_-)}{\mathcal{L}^2(\Omega^0_-)}$
over time.

\subsubsection{Benchmark problem II}

The physical parameters for the test case 2 in \cite[Table~I]{HysingTKPBGT09}
are given by
\begin{equation} \label{eq:Hysing2}
\rho_+ = 10^3\,,\quad \rho_- = 1\,,\quad \mu_+ = 10\,,\quad \mu_- = 0.1\,,\quad
\gamma = 1.96\,.
\end{equation}
We report on the quantitative results for this benchmark computation in
Table~\ref{tab:risingbubbleII}, where we have used
$J_\Gamma=128$ interface elements and $J_\Omega^0=35\,092$ initial
bulk elements. The final number of bulk elements, $J_\Omega^M$, for the schemes
$(\schemeAex)$, $(\schemeAim)$ and $(\schemeALE)$ was 10\,638, 10\,716 and
9\,972, respectively.
\begin{table*}
\center
\hspace*{-3.25cm}
\begin{tabular}{rrrr}
\hline
& ($\schemeAex)$ & $(\schemeAim)$ & $(\schemeALE)$ \\
\hline
$\strikes_{\min}$                & 0.5435 & 0.5473 & 0.5266 \\
$t_{\strikes = \strikes_{\min}}$ & 3.0000 & 3.0000 & 3.0000 \\
$V_{c,\max 1}$                   & 0.2503 & 0.2502 & 0.2502 \\
$t_{V_c = V_{c,\max 1}}$         & 0.7290 & 0.7290 & 0.7300 \\
$V_{c,\max 2}$                   & 0.2403 & 0.2403 & 0.2400 \\
$t_{V_c = V_{c,\max 2}}$         & 2.1070 & 2.1540 & 2.0670 \\
$z_c(t=3)$                       & 1.1383 & 1.1387 & 1.1385 \\
CPU[s]                           & 132\,498 & 236\,635 & 236\,253 \\
\hline
\end{tabular}
\hspace*{-3.25cm}
\caption{($\rho_+ = 10^3,\rho_- = 1,\mu_+ = 10,\mu_- =0.1,\gamma = 1.96$)
Benchmark problem II on $\Omega = (0,1) \times (0,2)$ over the time
interval $[0,3]$.}
\label{tab:risingbubbleII}
\end{table*}%
We observe that the results in Table~\ref{tab:risingbubbleII}
are in good agreement with the corresponding numbers from the finest
discretization run of group 3 in \cite{HysingTKPBGT09}, which are given by
0.5144, 3.0000, 0.2502, 0.7317, 0.2393, 2.0600 and 1.1376. Here we note that
there is little agreement on these results between the three groups in
\cite{HysingTKPBGT09}, but we believe the numbers of group 3 to be the most
reliable ones. In fact, very similar results are also obtained in
\cite[Tables~5 and 6]{fluidfbp}. Moreover, we note that the number of full
remeshings performed for this test case using the schemes $(\schemeAex)$,
$(\schemeAim)$ and $(\schemeALE)$ was 787, 788 and 808, respectively.
These large numbers are caused by the strong deformation of the rising bubble,
with the narrow tail structure leading to very elongated bulk elements. In
fact, until time $t=2$ only 6 remeshings are needed for each of the three
schemes.

It is interesting to note that for this benchmark, the methods in
\cite{HysingTKPBGT09} do not agree on whether the two tails of the rising
bubble experience pinch-off of some droplets. But we note that the most
accurate method there indicates that no pinch-off occurs, in agreement with our
results and the results in \cite{AlandV12,fluidfbp}. Of course, in a front
tracking method like ours topological changes need to be performed
heuristically when they occur, and this can be done, for example, as described
in \cite{BrochuB09,Sacconi15}.

\begin{figure}[htbp]
\centering
\subfloat[$t=10^{-1}$]{\includegraphics[width=.45\textwidth]
{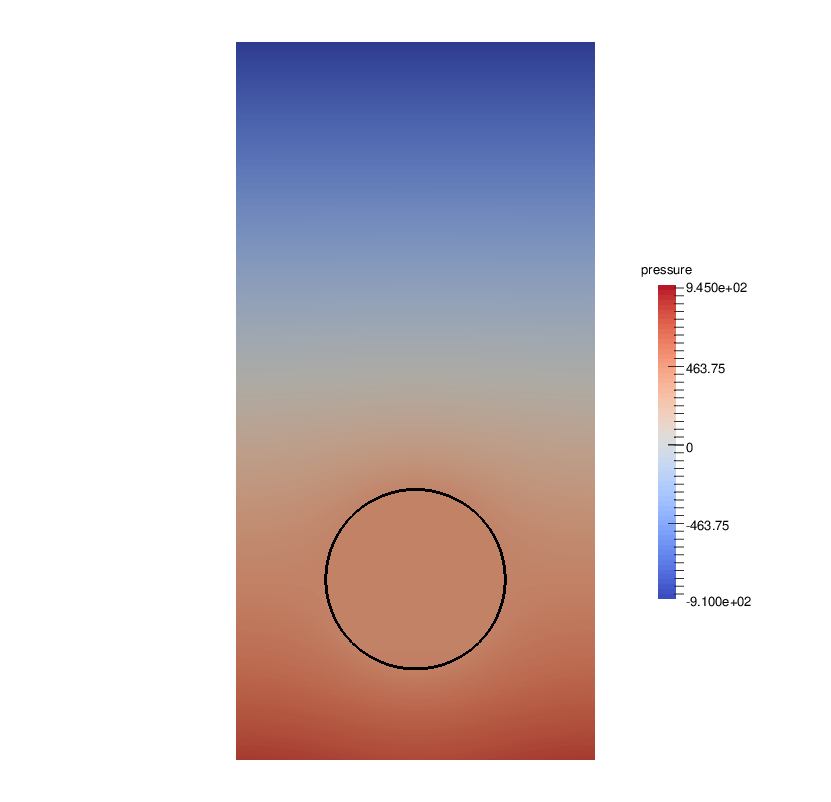}}
\subfloat[$t=1$]{\includegraphics[width=.45\textwidth]
{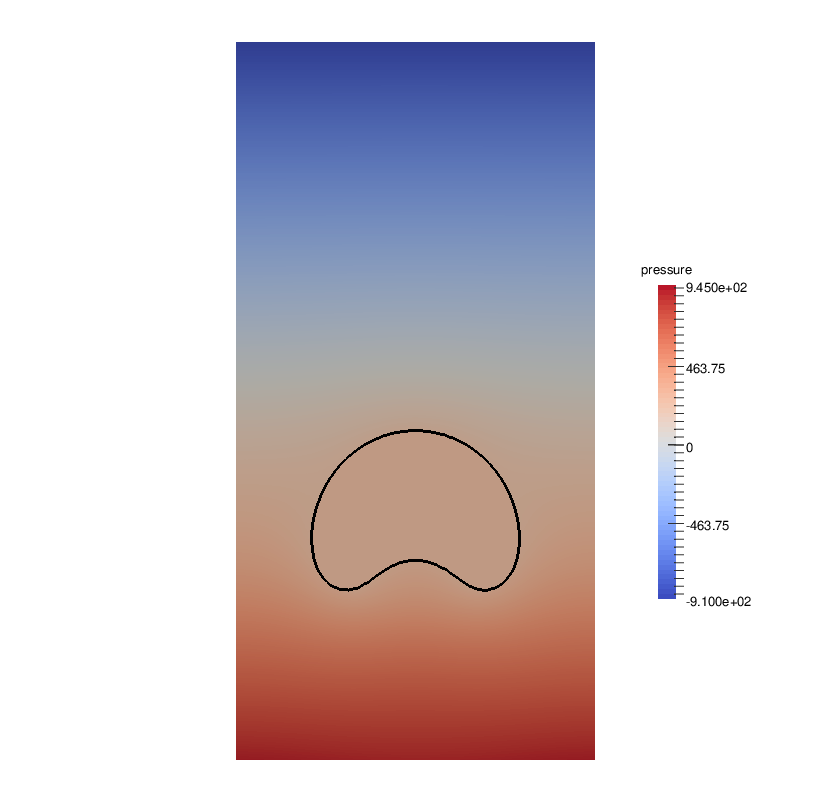}}\\
\subfloat[$t=2$]{\includegraphics[width=.45\textwidth]
{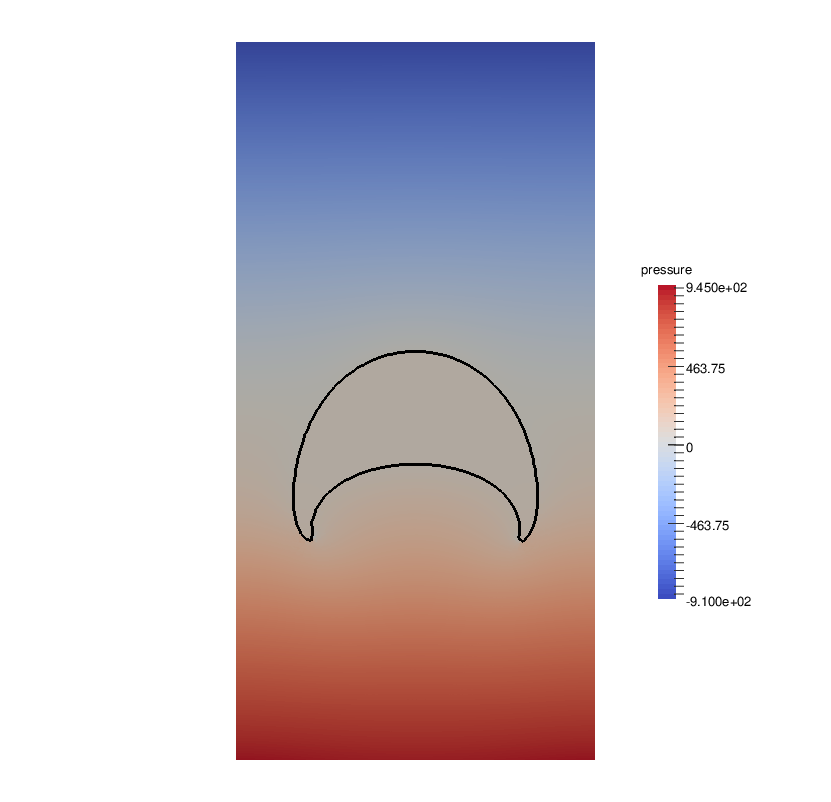}}
\subfloat[$t=3$]{\includegraphics[width=.45\textwidth]
{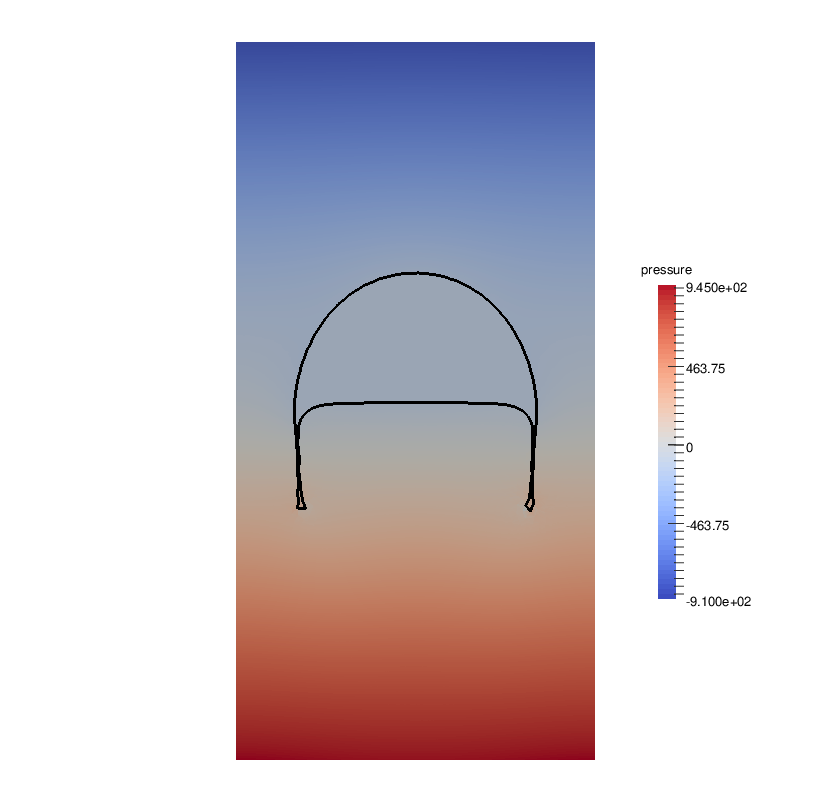}}
\caption{($\rho_+ = 10^3,\rho_- = 1,\mu_+ = 10,\mu_- =0.1,\gamma = 1.96$)
Pressure evolution for the benchmark problem II for the scheme $(\schemeAex)$,
with $J_\Gamma=128$ interface elements.}
\label{fig:risingbubbleIIpressure}
\end{figure}

\begin{figure}[htbp]
\centering
\subfloat[$t=10^{-1}$]{\includegraphics[width=.45\textwidth]
{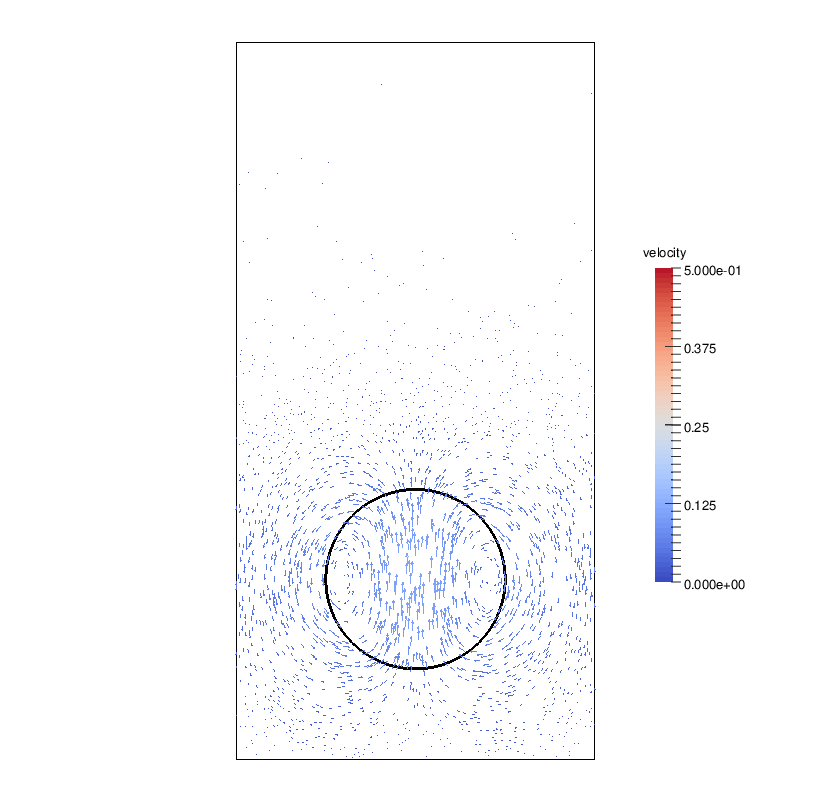}}
\subfloat[$t=1$]{\includegraphics[width=.45\textwidth]
{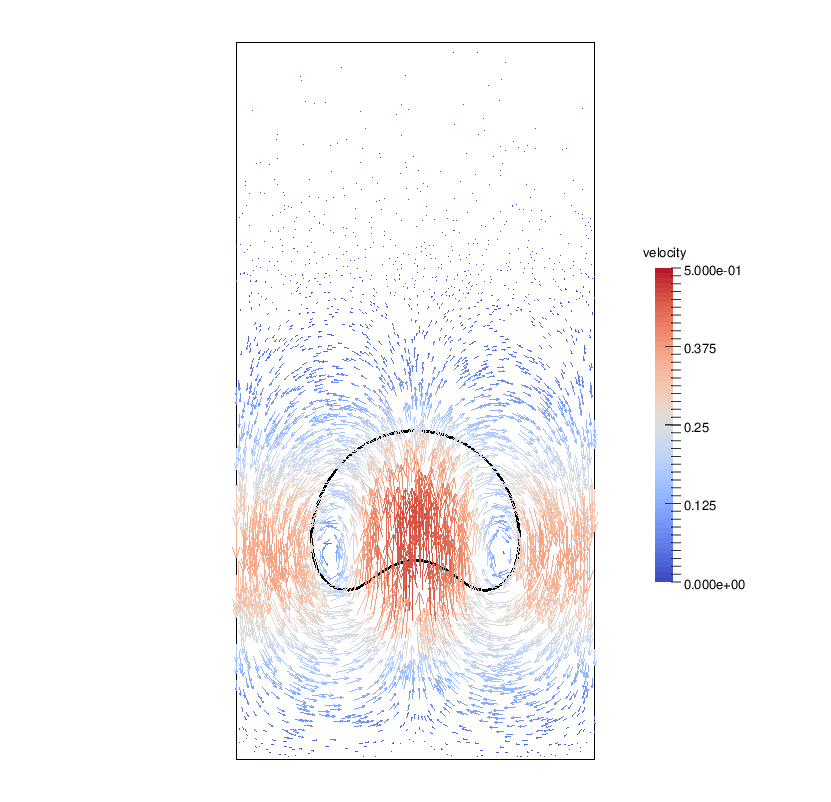}}\\
\subfloat[$t=2$]{\includegraphics[width=.45\textwidth]
{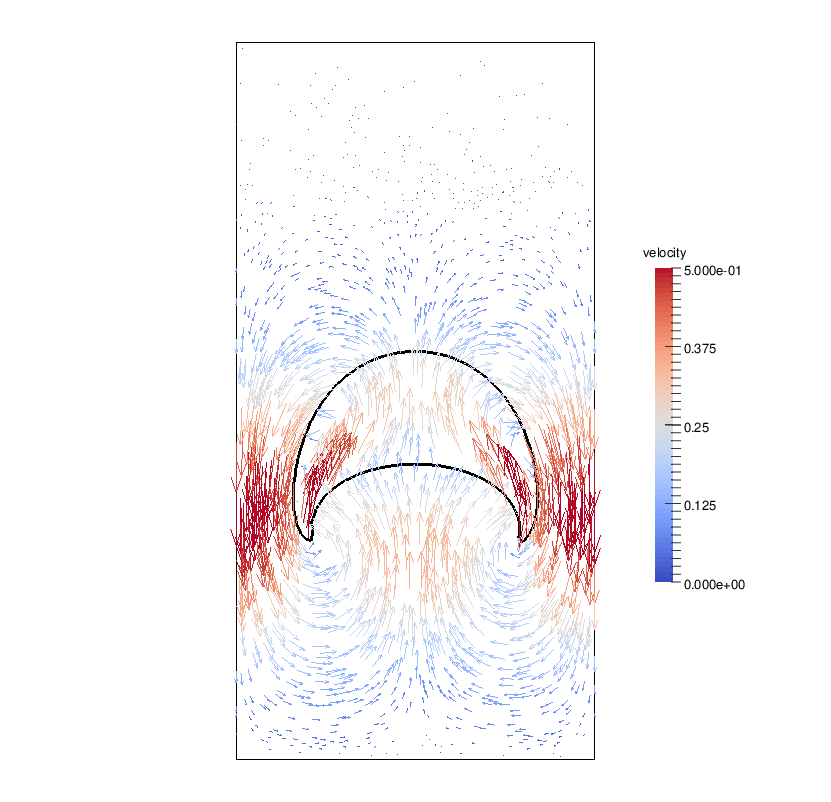}}
\subfloat[$t=3$]{\includegraphics[width=.45\textwidth]
{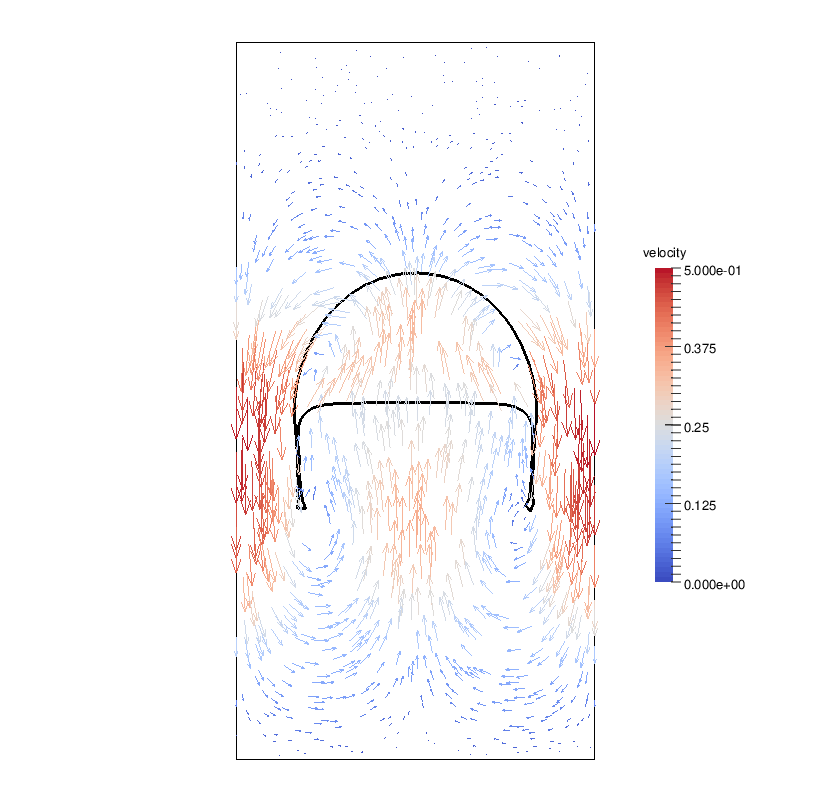}}
\caption{Velocity vector field for the simulation in
Figure~\ref{fig:risingbubbleIIpressure}.}
\label{fig:risingbubbleIIvelocity}
\end{figure}

\begin{figure}[htbp]
\centering
\includegraphics[width=.45\textwidth]
{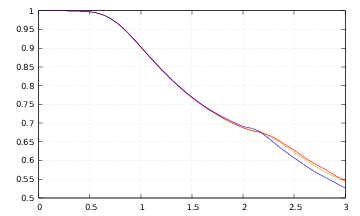}
\caption{A plot of the sphericity $\strikes$ over time for the simulation in
Figure~\ref{fig:risingbubbleIIpressure} for the schemes $(\schemeAex)$
(orange), $(\schemeAim)$ (red) and $(\schemeALE)$ (blue).}
\label{fig:risingbubbleIIsphericity}
\end{figure}

\begin{figure}[htbp]
\centering
\includegraphics[width=.45\textwidth]
{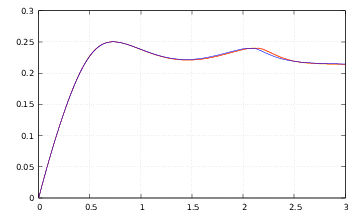}
\caption{A plot of the rising velocity $V_c$ over time for the simulation in
Figure~\ref{fig:risingbubbleIIpressure} for the schemes $(\schemeAex)$
(orange), $(\schemeAim)$ (red) and $(\schemeALE)$ (blue).}
\label{fig:risingbubbleIIrisingvelocity}
\end{figure}

\begin{figure}[htbp]
\centering
\includegraphics[width=.45\textwidth]
{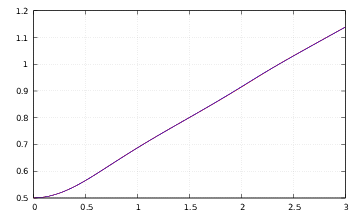}
\caption{A plot of the barycentre $z_c$ over time for the simulation in
Figure~\ref{fig:risingbubbleIIpressure} for the schemes $(\schemeAex)$
(orange), $(\schemeAim)$ (red) and $(\schemeALE)$ (blue).}
\label{fig:risingbubbleIIbarycenter}
\end{figure}

\begin{figure}[htbp]
\centering
\includegraphics[width=.45\textwidth]
{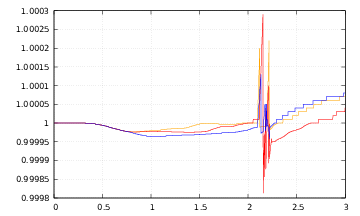}
\caption{A plot of the relative inner area
$\frac{\mathcal{L}^2(\Omega^m_-)}{\mathcal{L}^2(\Omega^0_-)}$ over time for the
simulation in Figure~\ref{fig:risingbubbleIIpressure} for the schemes
$(\schemeAex)$ (orange), $(\schemeAim)$ (red) and $(\schemeALE)$ (blue).}
\label{fig:risingbubbleIIinnervolume}
\end{figure}

In Figure~\ref{fig:risingbubbleIIpressure} we show the evolution of the discrete
pressures for a simulation with $J_\Gamma=128$ interface elements using the
scheme $(\schemeAex)$, while the velocities are visualized in
Figure~\ref{fig:risingbubbleIIvelocity}.
As before we show plots of the sphericity $\strikes$, the
rising velocity $V_c$, the barycentre $z_c$ and the relative
inner area $\frac{\mathcal{L}^2(\Omega^m_-)}{\mathcal{L}^2(\Omega^0_-)}$ over
time in Figures~\ref{fig:risingbubbleIIsphericity},
\ref{fig:risingbubbleIIrisingvelocity}, \ref{fig:risingbubbleIIbarycenter} and
\ref{fig:risingbubbleIIinnervolume}, respectively. In the former three
figures we note the good agreement with the plots shown in
\cite{HysingTKPBGT09}. As regards the plot of the relative inner area
in Figure~\ref{fig:risingbubbleIIinnervolume}, we observe oscillations in
the relative inner area between $t=2.1$ and $t=2.3$ which are caused by
the very thin filaments originating at the bottom of the rising bubble.

\section*{Conclusions}
We have introduced three Eulerian and one ALE scheme for the fitted finite
element approximation of two-phase Navier--Stokes flow. Computationally the
main difference between the Eulerian schemes and the ALE scheme is that
for the former a bulk mesh adjustment needs to be performed after every time
step, since the interface mesh in general moves, and so the fitted bulk mesh
needs to be adjusted. Here we employ a bulk mesh smoothing procedure based
on a linear elasticity problem. The mesh smoothing necessitates the
interpolation of the discrete velocity approximation onto the new bulk mesh.
It is often thought that the interpolation errors introduced by this procedure,
together with the additional computational effort for the calculation of the
interpolation, makes ALE schemes vastly superior to standard Eulerian
approximations. However, our numerical results do not bear this out.
In particular, for the physical and discretization parameters that we consider,
the Eulerian schemes $(\schemeAex)$ and $(\schemeAim)$ perform
as well as the ALE scheme $(\schemeALE)$ in most of our convergence tests
and benchmark computations. We conjecture that this good performance of
the Eulerian schemes is connected to the fact that all our schemes
maintain the interface mesh quality throughout, and so the
interface mesh, and the finite element information stored on it,
is never heuristically adjusted, not even during a full bulk remeshing.
Moreover, thanks to an efficient element-traversal
algorithm, the additional computational effort needed
for the velocity interpolation was only noticeable on nonconvex domains.


\def\soft#1{\leavevmode\setbox0=\hbox{h}\dimen7=\ht0\advance \dimen7
  by-1ex\relax\if t#1\relax\rlap{\raise.6\dimen7
  \hbox{\kern.3ex\char'47}}#1\relax\else\if T#1\relax
  \rlap{\raise.5\dimen7\hbox{\kern1.3ex\char'47}}#1\relax \else\if
  d#1\relax\rlap{\raise.5\dimen7\hbox{\kern.9ex \char'47}}#1\relax\else\if
  D#1\relax\rlap{\raise.5\dimen7 \hbox{\kern1.4ex\char'47}}#1\relax\else\if
  l#1\relax \rlap{\raise.5\dimen7\hbox{\kern.4ex\char'47}}#1\relax \else\if
  L#1\relax\rlap{\raise.5\dimen7\hbox{\kern.7ex
  \char'47}}#1\relax\else\message{accent \string\soft \space #1 not
  defined!}#1\relax\fi\fi\fi\fi\fi\fi}

\end{document}